\documentclass[11pt]{article}

\RequirePackage[OT1]{fontenc}
\RequirePackage{amsthm,amsmath}
\RequirePackage[numbers]{natbib}
\RequirePackage[colorlinks,citecolor=blue,urlcolor=blue]{hyperref}
\usepackage{amsfonts}
\usepackage{amssymb}
\usepackage{amsthm}
\usepackage{amsmath}
\usepackage{bigints}
\usepackage{algpseudocode}
\usepackage{bbm}
\usepackage{bm}
\usepackage{graphicx}
\usepackage{subcaption}
\usepackage{booktabs}
\usepackage{multirow}
\usepackage{longtable}
\usepackage[export]{adjustbox}
\usepackage{ragged2e}
\usepackage{threeparttable}

\numberwithin{equation}{section}
\theoremstyle{plain}
\newtheorem{theorem}{Theorem}
\newtheorem{definition}{Definition}
\newtheorem{Prop}{Proposition}
\newtheorem{corollary}{Corollary}
\newtheorem{example}{Example}
\newtheorem{lemma}{Lemma}

\newcommand{\Cov}{\mbox{Cov}}
\newcommand{\1}{\mathbbm{1}}
\newcommand{\R}{I\!\!R}
\newcommand{\N}{I\!\!N}

\newcommand{\E}{\mathbb{E}}

\usepackage{geometry}
\title{Covariance constraints for stochastic inverse problems of computer models}
\author{Nicolas Bousquet, M\'elanie Blaz\`ere, Thomas Cerbelaud}
\date{}

\begin{document}
\maketitle

\begin{abstract}
Stochastic inverse problems considered in this article consist of estimating the probability distributions of intrinsically random inputs of computer models. These estimations are based on observable outputs affected by model noise, and such problems are increasingly examined in parametric Bayesian contexts where the parameters of the targeted input distributions are affected by epistemic uncertainties. With the aim of improving the meaningfulness of solutions found by statistical algorithms --- in the sense that forward simulations based on such solutions must lead to relevant observables  --- we derive new prior constraints using the principles of global sensitivity analysis and information theory. Primarily formalized as constraints on covariances in Gaussian linear or linearizable situations, they reflect the idea that the solution should explain most of the observable uncertainty, while the model noise remains a secondary factor of this uncertainty. Simulated experiments highlight that, when injected into stochastic inversion algorithms, these constraints can indeed limit the influence of model noise on the result. They provide hope for future extensions in  more general frameworks, for example through the use of linear Gaussian mixtures.
\end{abstract}


\section{Introduction}

In the broadest sense, inverse problems consist of reconstructing a signal $X\in\R^p$ from  indirect observations $Y^*\in\R^q$ \cite{cavalier2011inverse}. \textcolor{black}{A special class of such problems is composed of stochastic inverse problems, where $X$ is considered to be intrinsically random, with an unknown or weakly known distribution $f_X$.   
The challenge here is to reconstruct $f_X$. This type of situation arises in many uncertainty quantification (UQ) problems in physics, where $X$ summarizes the behavior of non-measurable, uncertain variables such as physical correlations in thermo-hydraulics \cite{cocci2022comprehensive}, friction coefficients \cite{butler2015definition} or geometric and material properties \cite{choi2016stochastic,Sabbagh2015,Bingham2024}.  In this framework, the relationship between $Y^*$ and $X$ is often expressed by the following equation:
\begin{eqnarray}
Y^* & = & g(X) + \varepsilon, \label{sit1}
\end{eqnarray}
where $g:\R^p\to\R^q$ is a numerical model and $\varepsilon\in\R^q$ is a mixture of (epistemic) model error and (aleatoric) measurement error.  
Reconstructing $f_X$ is usually performed by statistical estimation from observations ${\bf y^*_n}=(y^*_1,\ldots,y^*_n)$ of $Y^*$, exploiting the missing data structure (the corresponding $x_i$) of the problem.  While nonparametric frameworks have been studied  \cite{Chafai2006,Bochkina2013,Vollmer2013,Dashti2013b}, most authors in statistical engineering select a parametrized form $f_X(x)=f(x|\theta)$ where $\theta$ lives in a  low-dimensional space $\Theta$ \cite{Barbillon2011,Fu2016,cocci2022comprehensive}. The problem therefore comes down to estimating $\theta$. This parametric choice meets practical needs: the number $n$ of available observations can be modest in real-life applications, and classical forms (as Gaussian or log-normal distributions \cite{cocci2022comprehensive}) allow to get identifiability conditions, at the cost of a linearization of $g$ \cite{Celeux2010}. For the same reason, estimation techniques exploiting solely the missing data structure, such as EM-type algorithms \cite{Celeux2010,Barbillon2011}, are often supplemented by Bayesian approaches involving Gibbs sampling \cite{Fu2015,Nagel2016,Barbillon2016,Dashti2013}, since they can incorporate additional prior knowledge on $\theta$. See \cite{Pereyra2015} for a recent survey of such algorithms and \cite{Delyon1999,Resmerita2007,Latuszynski2013} for details on their theoretical convergence properties.   
} 

\textcolor{black}{
In practice, the reconstruction of $f(x|\theta)$ can face two main challenges: 
\begin{enumerate} \item[(i)] Exploring a large $\Theta$ can cause computation to lose precision and place importance on implausible values $(\theta_k)_{k}$  \cite{Giryes2018}. This can result in underestimating the input uncertainty \cite{Kaipio2011,Dosso2012}, reflected in the features of $f(x|\theta)$, and more generally lead to irrelevant knowledge on $X$.
\item[(ii)] When $f(x|\theta)$ is reconstructed to be used for forward simulation-based predictive studies,  the simulations $Y:=g(X)$ can be far from the true observations. This can be due to Issue {(i)}, or/and because of a too strong model error within $\varepsilon$. The latter means that $g$ lacks of \emph{sufficient rigor for the context of forward use} \cite{Muguleta2018} and should be improved. But such improvements may conflict with the ease of use ($g$ can be deliberately chosen as a simplified or surrogate model, e.g.  \cite{chen2016model,keane2020robust}) and common validation methodologies are usually based on a specific, finite range of input/output particular values selected on considerations (e.g., limit cases) weakly related to the sampling mechanism $\{X\sim f(x|\theta)\}$ \cite{Oberkampf2006,Tedeschi2006,Oberkampf2010}. Ultimately, it is uneasy to differentiate the origin of such discrepancies between simulations and observations.    
\end{enumerate}
}
\textcolor{black}{To mitigate these difficulties, particular attention is generally given to use regularization tools, interpretable as prior constraints on $\theta$ \cite{Cavalier2011}. For example, if $\theta$
includes a covariance matrix, one might wish to impose sparsity. However, selecting finely such constraints is made difficult by the fact that $X$ is an unobserved multidimensional variable.  
Therefore one might prefer to propose constraints that indirectly apply to $X$ through considerations on observable variables. For physical reasons, if $g$ is sufficiently accurate, it can be reasonably assumed that the variability of $Y$ is close to that of $Y^*$. If this is not the case, the error $\varepsilon$ plays a dominant role in  (\ref{sit1}), threatening to limit the quality of the reconstruction of $f(x|\theta)$, as stressed in \cite{Fu2015}. 
Moreover, any forward use of $g(X)$ would be risky, potentially restricting or even prohibiting its application in a wide range of studies, such as optimization under uncertainty \cite{keane2020robust} or structural reliability analyses \cite{lemaire2013structural}. \\
How can we derive a prior constraint on $\theta$ from this rationale? A strong proximity between the variabilities of $Y$ and $Y^*$ necessarily requires that $X$ is a more influential variable than $\varepsilon$ for the variability of $Y^*$. 
This idea of ``influence" can be expressed using the concepts of global sensitivity analysis (GSA; \cite{da2021basics}), an essential UQ sub-domain fed by information theory \cite{Cover2012}: An uncertain input variable (or a group of variables) $X$ of $g$  is said to be influential on the output observable if it is a major factor in how the latter varies -- and equivalently, in its uncertainty.  Drawing inspiration from these concepts, Definition \ref{meaningful} is proposed to express this necessary condition. The term \emph{meaningful} is suggested to echo the idea that the simulation resulting from the inversion is necessarily close to the observations in a certain sense, and thus that the estimation of $f(x|\theta)$ should have practical relevance, conditional on the choice of $g$. }

\begin{definition}\label{meaningful}
The stochastic inverse problem (\ref{sit1}) is said \emph{meaningful} for $g$ if, given an error $\varepsilon\sim f_{\varepsilon}$, $f(x|\theta)$ carries more uncertainty over  $Y^*$ than $f_{\varepsilon}$.    
\end{definition}

\textcolor{black}{This definition remains fuzzy since ``carries more uncertainty" should be formalized. However, it can be adapted for different approaches to GSA, as brought to light later in the text. 
Building on this definition,  using tools from GSA seems appropriate to verify the result of an inversion before using it in simulations. It may also be interesting to leverage these tools to influence the inversion process itself, by generating the prior constraints evoked previously. First steps on this topic are investigated in this paper. \\
To our knowledge, this idea of building prior constraints based on GSA was not previously explored in the UQ field, although the need of constraining $\theta$ with respect to the features of an error $\varepsilon$ was  highlighted by \cite{Fu2015} about stochastic inversion problems involving Kriging-based surrogates. The reasons may be that historical approaches to inverse problems in the presence of multiple sources of uncertainty are difficult, still largely based on plug-in estimation \cite{Kennedy2001,HigdonKennedy,HigdonGattiker} to avoid identifiability problems. Beyond, stochastic inversion of highly-dimensional numerical models remains challenging to implement \cite{Damblin2015,Lee2019}.
}  

\textcolor{black}{
In the remainder of this paper, we focus on a first exploration of this idea by assuming, as many authors, that $X\sim {\cal{N}}_p(\mu,\Gamma)$ with $\theta:=(\mu,\Sigma)\in\R^p \times S^{++}_p(\R)$ where $S^{++}_p(\R)$ is the space of symmetric, positive-definite matrices in $\R^p$, and that $\varepsilon\sim {\cal{N}}_q(0,\Sigma)$ with $\Sigma\in S^{++}_q(\R)$. Besides, we consider first linear then linearizable inversion problems, motivated by their prominent role in this field.   From Definition \ref{meaningful} and GSA principles, we derive several possible formalizations leading mainly to explicit inequalities between the features of the input covariance $\Gamma$ and those of  $(g,\Sigma)$. Most of these inequalities  are  tested as prior constraints on $\Gamma$ within  stochastic inversion algorithms. These experiments are conducted in a Bayesian setting, as it appears as the most natural framework. Experimental results highlight that such constraints can help to improve the statistical estimation of $f(x|\theta)$ when this distribution is used for predictive analysis. Beyond these results, this technical work leads to a new characterization of Fisher information for Gaussian models, allowing us to understand how this information is degraded in the event of model error. It also introduces a link between the notion of influence proposed by GSA and that of signal-to-noise ratio, itself related to the theory of well-posed inversion problems in Hadamard's sense  \cite{kabanikhin2008definitions}.  
}

More precisely, this article is structured as follows. \textcolor{black}{Main notations and concepts are introduced in Section \ref{notations}.} In a preliminary approach, Section \ref{linear-1D} studies the case of a linear model $g$ and introduces two simple constraint rules based on well-known GSA tools, Sobol's and entropy-based indices, which provide similar results. 
\textcolor{black}{A short study is conducted in view of extending these results to locally linearizable models (requiring $g$ being locally differentiable), and illustrates the difficulty of working with such approximations.} Then Section \ref{section:well-posedness-FI} provides and studies a general constraining rule, based on Fisher information, in broader settings when $g$ can be approximated by a \emph{globally} linear model in a variational sense. Based on a real case-study previously investigated in the UQ literature, numerical experiments previously evoked are conducted on Section \ref{numeric}. 
A discussion section ends this article, opening avenues for future research. Proofs of the stated results are provided in a dedicated Appendix.

\section{Formalism and main notations}
\label{notations}

\textcolor{black}{
We introduce here the formalism necessary for the mathematical developments in the following sections. For $p\in\N^*$, $p<\infty$, denote ${\cal{B}}(\R^p)$ the Borel $\sigma-$algebra that endows the metric space $(\R^p,d)$ where $d$ is induced by the Euclidian norm  (generically noted $\|.\|$ along the article, and sometimes $\|.\|_p$ to reflect the appropriate dimension and avoid ambiguity).}

Let $(\R^p,{\cal{B}}(\R^p),\mathbb{P})$ and $(\R^q,{\cal{B}}(\R^q),P)$ be two probability spaces, $E=\otimes_{i=1}^p E_i$ be endowed with a product $\sigma-$algebra and $(X,\varepsilon)$ be two  independent random variables from $(\R^p, \R^q)$ to $(E,\R^q)$ in $L^2$. Finally, let $g:(E,\R^p)\to(\R^q,{\cal{B}}(\R^q))$ be in $L^2$ and
 \begin{eqnarray*}
 Y & = & g(X), \\
 Y^* & = & Y + \varepsilon
 \end{eqnarray*}
be two random variables defined on  $(\R^q,{\cal{B}}(\R^q),P)$.   
To use classical results in GSA, we also assume that $\E[Y^2]<\infty$. 

\textcolor{black}{On the basis of this formalism, and within the framework precised in Introduction, we consider the following stochastic inversion problem:
\begin{eqnarray}
& \text{Given a sample of observations ${\bf y^*_n}=(y^*_1,\ldots,y^*_n)$ of $Y^*$,} \nonumber \\
& \ \ \ \ \ \ \ \text{find $\theta=(\mu,\Gamma)\in \R^p\times S^{++}_p(\R)$ such that, for $i\in\{1,\ldots,n\}$:} \nonumber \\
& \left\{
\begin{array}{lll}
 x_i & \sim & {\cal{N}}_p(\mu,\Gamma), \\
 y^*_i & = &  g(x_i) + \varepsilon_i
 \end{array}
 \right. \label{inversion.problem} \\
& \ \ \ \  \text{where $\varepsilon_1,\ldots,\varepsilon_n \overset{\text{\tiny iid}}{\sim} {\cal{N}}_q(0,\Sigma)$ and $\Sigma\in S^{++}_q(\R)$ is known. }
\nonumber
\end{eqnarray}
}

In the remainder of this article, 
the symbol $|B|$ will be used for denoting the determinant of matrix $B$. If $B$ is a diagonalizable $p\times p$ matrix, denote $\lambda^{B}_1,\ldots,\lambda^{B}_p$ the set of its eigenvalues. Furthermore, we use the notation $I_q$ for the unit matrix in $\R^{q,q}$. Finally, given a set $\{x_1,\ldots,x_n\}$ of real values, denote 
 $x_{\inf}=\min_{1\leq i \leq n} x_i$ and $x_{\sup}=\max_{1\leq i \leq n} x_i$.

\section{A preliminary: Two intuitive notions of meaningful solutions}
\label{linear-1D}

\subsection{Linear models}\label{linear.A}

A first approach to formalizing Definition \ref{meaningful} can be proposed using the case of the simple linear model
\begin{eqnarray} 
g(x)=a \cdot x \ \ \ \ \text{with $a\in\R^{q,p}$.} 
\label{eq:model_toy}
\end{eqnarray} 

 \noindent In GSA, covariance and entropy play important roles in the comparison of uncertainties modeled by probability distributions. Related to output covariance and recalled in Appendix \ref{sobol}, the first-order grouped Sobol' indices seem appropriate candidates to express 
the fact that most of the uncertainty on $Y^*$ must be explained by the uncertainty on $X$, accounting for possible correlations among the components of $X$. At the price of assuming independence between $X$ and $\varepsilon$, Definition \ref{soboldef} is a possible formalization of Definition \ref{meaningful}.

\begin{definition}\label{soboldef}
 Let $(S_{X},S_{\varepsilon})$  be the  \textcolor{black}{first-order grouped Sobol indices} quantifying the uncertainty on $Y^*$ explained by $X$ and $\varepsilon$, respectively. Provided $X$ and $\varepsilon$ are independent, the stochastic inverse problem (\ref{inversion.problem}-\ref{eq:model_toy}) is said to be meaningful in Sobol' sense if 
\begin{eqnarray} 
 S_X  & > &  S_{\varepsilon}. \label{sobol1}
 \end{eqnarray} 
 \end{definition} 

The entropy-based indices proposed by \cite{KH01} and studied by \cite{Liu2006,Auder2009} offer an alternative view. Denoting $\mathcal{E}(Y^*)$ the Shannon-Kullback entropy of $Y^*$, this approach stands on the notion of {conditional entropy} $\E_X[\mathcal{E}(Y^*)|X]$, that measures the average loss of information on $Y^*$ when the behavior of $X$ is known. If $\varepsilon$ is of small influence in the variability of $Y^*$, then it can be expected that  $\E_X\left[\mathcal{E}(Y^*)|X\right]<\E_{\varepsilon}\left[\mathcal{E}(Y^*)|\varepsilon\right]$. This leads to the following alternative formalization of Definition \ref{meaningful}. \\

\begin{definition}\label{entropy} 
The stochastic inverse problem (\ref{inversion.problem}-\ref{eq:model_toy}) is said to be meaningful in the entropic sense if 
\begin{eqnarray}
\E_{\varepsilon}\left[\mathcal{E}(Y^*)|\varepsilon\right] & > & \E_X\left[\mathcal{E}(Y^*)|X\right]. \label{entropy1}
\end{eqnarray} 
 \end{definition}

Next proposition straightforwardly provides prior constraints on $\Gamma$ from Definitions \ref{soboldef} and \ref{entropy}.

\begin{Prop}\label{prop1_bis}
The  stochastic inverse problem   (\ref{inversion.problem}-\ref{eq:model_toy}) is meaningful in Sobol' sense if and only if
\begin{eqnarray}
{\rm{Tr}}( a\Gamma a^T)   >  {\rm{Tr}}(\Sigma). \label{eq:model2bis1}
\end{eqnarray} 
It is meaningful in entropic sense if and only if
\begin{eqnarray}
\left|a\Gamma a^T\right|   >  \left|\Sigma\right|. \label{eq:model2ter1}
\end{eqnarray} 
which imposes that $a\Gamma a^T$ be invertible.
\end{Prop}

\noindent When the symmetric matrix $a\Gamma a^T$ is close to the identity, ``the determinant behaves like the trace" \cite{Tao2013}: for any bounded matrix $A\in\R^{q, q}$ and infinitesimal $\varepsilon'\in\R$, denoting $I_q$ the unit matrix of rank $q$, 
\begin{eqnarray*}
| I_q + \varepsilon' A|  & = & 1 + \varepsilon' {\rm{Tr}}(A) + O(\varepsilon'^2).
\end{eqnarray*}
In such cases (\ref{eq:model2bis1}) and (\ref{eq:model2ter1}) reflect a similar constraint on $\Gamma$. 

Besides, an immediate corollary of Proposition \ref{prop1_bis} is that both results coincide when the output is univariate. \\

\begin{corollary}
When $q=1$, the stochastic inverse problem (\ref{inversion.problem}-\ref{eq:model_toy}) is meaningful in Sobol' and entropic senses if and only if
\begin{eqnarray}
a^T\Gamma a &  > & \sigma^2. \label{eq:model2}
\end{eqnarray}
\end{corollary}

In a Bayesian stochastic inversion framework, an interesting result can be derived, under the additional assumptions that $q=p$  (then  (\ref{eq:model2ter1}) becomes $\left|a\Gamma a^T\right|=|a|^2 \left|\Gamma\right|$), $a$ be invertible, and that there exists $\beta<\infty$ such that $|\Gamma|<\beta$. \textcolor{black}{Such a constraint is likely to appear in practical applications. For instance, consider that for all $i\in\{1,\ldots,p\}$, $\rm{Var}(X_i)\leq \beta_j$. Then, from Hadamard's inequality, $\beta=\prod_{i=1}^p\beta_i$.  The following proposition deals with the possibility of handling an integrable Jeffreys measure as a prior for $\Gamma$, and thus providing a practical objective framework for the Bayesian estimation of this covariance matrix. For any parametric density model $f(x|\theta)$ with well-defined Fisher information $I(\theta)$ (see Appendix \ref{fisher} for details), Jeffreys' measure is considered as one of the reference non-informative prior measures on $\theta$ in Bayesian estimation, with rich theoretical properties \cite{Clarke1994}. Especially, it is invariant through any bijective reparameterisation of $\theta$.  Jeffreys' prior is defined by
\begin{eqnarray*}
\pi^J(\theta) & \propto & \sqrt{|I(\theta)|}. 
\end{eqnarray*}
However, it is often an improper (non-integrable) measure, which can lead to posterior inconsistencies and paradoxes \cite{Kass1996}.This is especially the case where no lower bound other than 0 is available for  $|\Gamma|$. 
Next result (Proposition \ref{jeffreys}) is therefore likely to be useful in contexts of objective Bayesian modeling and selection \cite{Berger2001}.} 

\begin{Prop}\label{jeffreys}
Assume  $q=p$ and $|a|>0$. Let 
\begin{eqnarray*}
\Omega & = & \left\{\Gamma \in S^{++}_p(\R) \ ; \ 0<\frac{|\Sigma|}{|a|^2} < |\Gamma| < \beta\right\}
\end{eqnarray*}
be the supporting set for the estimation of covariance matrix $\Gamma$ under Definition \ref{entropy}. Let $A$ denote a compact set of $\R^p$. Then Jeffrey's prior measure for  $(\mu,\Gamma)$ in (\ref{inversion.problem}-\ref{eq:model_toy}), restricted to  $A\times\Omega$, is a proper (integrable) probability measure.
\end{Prop}

\subsection{Locally linearizable models}\label{locally.linear}

\noindent  While  linear models  are generally used in a pure explanatory regression context,  one should consider rather, to address more broader situations and examine how to adapt the previous covariance constraints, that $g$ is  linearizable in the neighborhood of  a given point of interest $x_0$.
Such a \emph{local} linearization of $g$ was defended by \cite{Barbillon2011}, among others,  to avoid prohibitive computational costs in simulation-based inferential algorithms.  More generally, local linearization is an usual tool for dealing with inverse problems \cite{Stefanov2009}. \\

\textcolor{black}{Local linearization can be studied using the multivariate first-order Taylor expansion of the vector-valued function $g$ 
 under the assumption that $g$ is Fr\'echet differentiable in a neighborhood $V_{x_0}$ of a given $x_0\in\R^p$. More precisely, denote
 $g^{(i)}$ the $i$th coordinate of $g$ and let $Dg_{x_0}=(\frac{\partial g^{(i)}}{\partial x_j} (x_0))_{i,j}$ be the $(q,p)$ Jacobian matrix of $g$.
 \begin{description}
\item[H0]: $g:\R^p\to \R^q$ is such that $\forall i\in\{1,\ldots,q\}$, $\forall j \in\{1,\ldots,p\}$,   $\frac{\partial g^{(i)}}{\partial x_{ij}}$ exists almost everywhere  and $\E\left|\frac{\partial g^{(i)}}{\partial x_j}(X) \right|<\infty$ on $V_{x_0}$. 
\end{description}
 Under H0, for any $x\in V_{x_0}$, one may write $g(x)  =  \tilde{g}_{x_0}(x) + h_{x_0}(x)$ with
\begin{eqnarray}
 \tilde{g}_{x_0}(x) & = & g(x_0) + Dg_{x_0}(x-x_0) \ = \ g(x_0) + \left(\sum_{j=1}^{p}\dfrac{\partial g^{(i)}}{\partial x_j}(x_0)\left( x_j-x_{0_j} \right)\right)_{1\leq i \leq q} \label{first.order.approx}
\end{eqnarray}
and $h_{x_0}:\R^p\to \R^q$ is a remainder such that  $h_{x_0}(x)=o\left( \Vert x - x_0\Vert_p\right)$. }\\

\textcolor{black}{
Obtaining a constraint on the covariance of $g(X)$ from that of $\tilde{g}_{x_0}(X)$ requires modifying the rules (\ref{eq:model2bis1}-\ref{eq:model2ter1}) by involving quantities related to the geometry of $g(X)$ in $V_{x_0}$. In next Proposition we focus on Definition \ref{soboldef} and consider the case of an univariate output $(q=1)$. 
}

\begin{Prop}\label{prop.gruff}
Assume H0 and $q=1$.  
Then a sufficient condition for (\ref{sobol1}) is 
 \begin{eqnarray}
{\rm{Tr}}\left[D_{g_{x_0}}\Gamma\left(2\E\left[D_g(X)\right] - D^T_{g_{x_0}}\right) \right] & \geq &  {\rm{Tr}}\left[\Sigma\right]. \label{res.propa}\label{constraint.adi}
\end{eqnarray}
\end{Prop}

\vspace{0.5cm}

\noindent
Note that when $g(X)=\tilde{g}_{x_0}(X)=a X + b$ with $a\in\R^{q,p}$ and $b\in\R^q$, then $D^T_{g_{x_0}}=a^T$ and  (\ref{res.propa}) is exactly the necessary and sufficient condition (\ref{eq:model2bis1}).
Next example illustrates the fact that (\ref{constraint.adi}) becomes actually a condition on the full parameter vector $(\mu,\Gamma)$, because of  linearization. 

\begin{example}\label{local.linear}
Assume $g(x)=1-\exp(-\sum_{i=1}^p x_i)$. Then, around $x_0=(0,\ldots,0)$, consider $\tilde{g}_{x_0}(x)=\sum_{i=1}^p x_i$. Denoting $\bar{\mu}=p^{-1}\sum_{i=1}^p \mu_i$ and $\gamma^2=\sum_{i,j} \Gamma_{i,j}$, then $\tilde{g}_{x_0}(X)\sim{\cal{N}}(p\bar{\mu},\gamma^2)$ and $1-g(X)\sim{\cal{LN}}(-p\mu,\gamma^2)$. With $\E[D_g(X)]=\exp(-p\bar{\mu}+\gamma^2/2)u^T$ and $D_{g_{x_0}}=u$ where $u$ is the vector in $\R^p$ such that $u_i=1$ for all $i=1,\ldots,p$, then (\ref{constraint.adi}) becomes
\begin{eqnarray*}
\gamma^2\left\{2\exp\left(\gamma^2/2-p\bar{\mu}\right)-1\right\} & \geq & {\rm{Tr}}\left[\Sigma\right].
\end{eqnarray*}
\end{example}

\noindent Extending the  result of Proposition \ref{prop.gruff} to $q>1$ seems uneasy, although some refinements could probably be derived using trace inequalities (see the proof). But this is enough to motivate the following comment. Since we only dispose of noisy observations $y^*_1,\ldots,y^*_n$ of $g(X)+\varepsilon$, estimating consistently 
 $\E[D_g(X)]$ or similar functions involving a mean of the Jacobian prior to model inversion appears challenging, if not impossible when $q\neq p$ and $g$ is not an invertible function. Rather,   $\E[D_g(X)]$ is a function of $\theta=(\mu,\Gamma)$:
 \begin{eqnarray*}
 \E\left[D_g(X)\right] & = & \int_{\R^p} D_g(x) f_X(x|\theta) \ dx 
 \end{eqnarray*}
 and (\ref{constraint.adi}) should be inserted as a constraint on $\Theta$ within the stochastic inversion algorithms evoked hereinbefore. For each candidate $\theta_i$ in such routines, the integral above has to be computed numerically (e.g., by Monte Carlo batch sampling or gradient descent techniques). Even in the case where $g$ can be automatically differentiated, it is likely that the rejection rate of candidates $\theta_i,\theta_{i+1},\ldots$ will make the computational cost prohibitive. For this reason, variational linear approximation will be studied later in the article, instead of local linearization, as it involves the first and second-order moments of the distribution of $g(X)$, which can be consistently estimated from the avai\-lable observations.

\section{Fisher inversion constraints}\label{section:well-posedness-FI}

Definition \ref{soboldef} is based on first-order hierarchization of conditional variances. This appears somewhat limited to 
reflect how input uncertainty from $X$ or $\varepsilon$ is transmitted to the observed output $Y^*$.  
\textcolor{black}{Since information can be viewed as a reduction of uncertainty  \cite{Cole1993}, the main concepts of classical information theory \cite{Cover2012} appear useful in providing a broader framework for characterizing the effects of uncertainty.} 
Because the inverse problem is to estimate $\theta=(\mu,\Gamma)$, using a truly parametric measure of information seems more appropriate. The most usual measure of information is Fisher information \cite{Cover2012}. See Appendix \ref{fisher} for a technical reminder and details about its interpretation.

\subsection{Formalisation}\label{work.hypo.0}

 Under the conditions of existence (cf. Appendix \ref{fisher}), denote by $I_{g(X)}(\theta)\in S^{+}_q(\R)$ and $I_{Y^*}(\theta)\in S^{+}_q(\R)$ the Fisher information matrices provided respectively by $g(X)$ and $Y^*$ about $\theta$. Since the impact of $\varepsilon$ is to degrade information, \textcolor{black}{it can be expected that $I_{g(X)}(\theta)$ be greater than $I_{Y^*}(\theta)$, according to some order relation $\succeq $ defined for symmetric matrices. The Loewner order
  appears as the most natural candidate for $\succeq$, as it is the partial order induced by the translations over the  convex cone of positive semi-definite matrices. 
 It is therefore expected that
\begin{eqnarray}
\label{eq:fisher-condition}
 I_{g(X)}(\theta) & \succeq & I_{Y^*}(\theta)  \label{result1}
\end{eqnarray}
where $A \succeq B$, for two squared matrices $A$ and $B$, means that $A-B$ is a positive semidefinite matrix.} \\

Stating that most of information on $\theta$ in $Y^*$ is transmitted from $g(X)$ (or equivalenty from $X$)  implies that the difference between $I_{g(X)}(\theta)$ and $I_{Y^*}(\theta)$, which is a measure of the information loss because of the noise $\varepsilon$, should not be greater than a fraction $(1-1/c) I_{g(X)}(\theta)$ where $c>1$. Therefore we propose the following definition. 

\begin{definition}\label{fisher.wp}
For $c>1$, the stochastic inverse problem is meaningful in Fisher' sense if
\begin{eqnarray}
\label{eq:fisher-condition2}
I_{g(X)}(\theta) \succeq I_{Y^*}(\theta) \succeq \dfrac{1}{c}I_{g(X)}(\theta).
\end{eqnarray}
\end{definition}

An intuitive value for $c$ is 2, but further developments will show that $c$ can be related to important concepts in matrix inversion and signal processing, as linear system conditioning and signal over noise ratio. This can help to assess its value in practice. \\

 The following sections consider several usual situations encountered in stochastic inversion and provide necessary conditions (NC), sufficient conditions (SC) and necessary and sufficient conditions (NSC) on the features of the covariance $\Gamma$ of $X$ such that  (\ref{eq:fisher-condition2}) be verified. In practice, NC and NSC should be preferred to SC, as the latter correspond to prior constraints on $\Gamma$ that threaten to overestimate how $X$ varies.

\subsection{Gaussian linear models}\label{work.hypo}

We assume a linear $g:x\mapsto Hx$, with $H\in \mathbb{R}^{q , p}$ of full rank.  
By abuse of notation, we denote in next equations $\partial / \partial \mu$ and $\partial / \partial \Gamma$ the differentiations of functions of $\theta=(\mu,\Gamma)$ with respect to the free components of $\theta$. 
  We remind that if $A$ is real square matrix such that $A\succeq 0$, then $\lambda^A_{\sup} = \| A\|_2$ is the Euclidian norm. 

\subsubsection{Fisher matrix characterization}\label{fisher.carac}

 In order to provide useful conditions for (\ref{eq:fisher-condition2}), we first provide a characterization of the Fisher matrix for the Gaussian model
 \begin{eqnarray}
 Y_{\alpha} & = HX + \alpha \varepsilon \label{gen.model} \ \ \ \text{with $\alpha\in\{0,1\}$,}
 \end{eqnarray}
that encompasses the cases $Y_{0}=Y$ and $Y_1=Y^*$. Assuming $\Sigma \in S^{++}_q(\R)$, we denote 
 \begin{eqnarray*}
 \Psi_{\alpha}(\Gamma) & = & \Sigma^{-\alpha/2} H \Gamma H^T \Sigma^{-\alpha/2}
 \end{eqnarray*}
 (with the convention $ \Sigma^{0}=I_q$), which can be seen as a multivariate generalization of the {\it  signal over noise ratio} (MSNR) characterizing (\ref{gen.model}) ; see \cite{Czanner2008} for details on this concept widely used in signal processing and the characterization of channel features. 
 
\begin{theorem}
\label{theo:fisher-info-gen}
Denote $I_{{Y}_{\alpha}}(\theta)$ the Fisher information matrix for (\ref{gen.model}), assuming $\theta=(\mu,\Gamma)\in\R^p\times S^{++}_p(\R)$. Let ${\partial^{2}}/{\partial \Gamma^2}$ indicate the double differentiation with respect to the  free components of $\Gamma$. Then, assuming that $ \alpha I_q + \Psi_{\alpha}(\Gamma)$ is inversible,
\begin{eqnarray*}
I_{{Y}_{\alpha}}(\theta) & = & \left(
\begin{array}{ll}
I_{{Y}_{\alpha}}(\mu) & {\bf 0} \\
{\bf 0} & I_{{Y}_{\alpha}}(\Gamma)
\end{array}\right)
\end{eqnarray*}
where 
\begin{eqnarray}
I_{{Y}_{\alpha}}(\mu) & = & H^T\Sigma^{-\alpha/2} \left(\alpha I_q + \Psi_{\alpha}(\Gamma)\right)^{-1} \Sigma^{-\alpha/2}  H, \nonumber \\
I_{{Y}_{\alpha}}(\Gamma) & = & -\dfrac{1}{2}\sum\limits_{i=1}^q \dfrac{\partial^{2}}{\partial \Gamma^2} \log\left(\alpha + \lambda^{\Psi_{\alpha}(\Gamma)}_i\right).\label{decomp.Igamma}
\end{eqnarray}
\end{theorem}

\noindent Next proposition and its corollaries provide more refined results for the Fisher submatrices  in the frequent situation where $\Gamma$ is chosen to be diagonal. It highlights the role of the eigenvalues of the MSNR in the information matrix. It requires to introduce the semidefinite matrices $H_k=[h_{i,k}]^T_{1\leq i \leq q}$ and
 \begin{eqnarray}
\tilde{A}_{k,\alpha}  & = & \Sigma^{-\alpha/2} H_k H^T_k \Sigma^{-\alpha/2} \label{tilde.A.k}
 \end{eqnarray}
 that play a key role in the characterization of this Fisher matrix.

\begin{Prop}\label{going.further}
Assume $\Gamma=diag(\tau^2_1,\ldots,\tau^2_p)$ and $\theta=(\mu,\Gamma)\in\R^p\times S^{++}_p(\R)$. Denote $V^{\Gamma}_{\alpha}=(v_{\alpha,1},\ldots,v_{\alpha,q})$ an orthonormal basis of $(q,1)-$eigenvectors  for $\Psi_{\alpha}(\Gamma)$, and define
\begin{eqnarray*}
\nu^{\alpha}_{i,k}  =  v^T_{\alpha,i}\tilde{A}_{k,\alpha} v_{\alpha,i} & \text{and} & \eta^{\alpha}_{i,k,l} = v^T_{\alpha,i} \tilde{A}_{l,\alpha} \Psi^+_i(\Gamma)  \tilde{A}_{k,\alpha} v_{\alpha,i}
\end{eqnarray*}
where $\Psi_i(\Gamma) =  \lambda^{\Psi_{\alpha}(\Gamma)}_i I_q- \Psi_{\alpha}(\Gamma)$ and $A^+$ is the generalized Moore-Penrose pseudo-inverse of $A$.
Then
\begin{eqnarray}
I_{{Y}_{\alpha}}(\mu)       & = &  H^T \Sigma^{-\alpha/2}  \left(\alpha I_q + \sum\limits_{k=1}^p \tau^2_k \tilde{A}_k\right)^{-1} \Sigma^{-\alpha/2}  H, \label{fisher.term.simp.1} 
\end{eqnarray}
and $\forall (k,l)\in\{1,\ldots,p\}^2$, the $(k,l)-$term of matrix $I_{{Y}_{\alpha}}(\Gamma)$ is
\begin{eqnarray}
I_{{Y}_{\alpha}}(\Gamma)_{k,l} & = & \dfrac{1}{2}\sum\limits_{i=1}^q \left\{ \frac{\nu^{\alpha}_{i,k}\nu^{\alpha}_{i,l}}{\left(\alpha + \lambda_i^{\Psi_{\alpha}(\Gamma)}\right)^2} - \frac{\eta^{\alpha}_{i,k,l} + \eta^{\alpha}_{i,l,k}}{\alpha + \lambda_i^{\Psi_{\alpha}(\Gamma)}}\right\}. \label{fisher.term.simp.2} 
\end{eqnarray}
\end{Prop}

\noindent Next corollary assumes that each outer product $H_k H^T_l$ is symmetric, for all couples $(k,l)\in\{1,\ldots,p\}$. It is  verified when $q=1$ (univariate output), when $p=1$ (univariate input, since $H_k=H_l$)  or more generally when the rank of $H$ is $1$. 
In the special case where $p=q$ and $H$ is a diagonal matrix, then it is easy to see that $H_k H^T_l$ is symmetric too. 

\begin{corollary}\label{diagonalite.aa}
Assume $\Gamma=diag(\tau^2_1,\ldots,\tau^2_p)$ and $\theta=(\mu,\Gamma)\in\R^p\times S^{++}_p(\R)$. Assume furthermore that $\Sigma$ is diagonal and $H_k H^T_l$ is symmetric for all couples $(k,l)\in\{1,\ldots,p\}^2$. Denote $\Lambda_{{\alpha}}=(\| \Sigma^{-\alpha/2} H_1\|^2_2,\ldots, \| \Sigma^{-\alpha/2} H_p\|^2_2)^T$ and  $\tilde{\Lambda}_{{\alpha}}=diag \Lambda_{{\alpha}}$. 
Then
\begin{eqnarray}
I_{{Y}_{\alpha}}(\Gamma) & = & \frac{\Lambda_{{\alpha}} \Lambda^T_{{\alpha}}}{2\left(\alpha + {\rm{Tr}}\left(\tilde{\Lambda}_{\alpha} \Gamma\right) \right)^2}. \label{I.mu.cov}
\end{eqnarray}
\end{corollary}

 Without the assumption of a diagonal $\Sigma$, Corollary \ref{diagonalite.1} states a last useful result derived straightforwardly from Theorem \ref{theo:fisher-info-gen}, assuming $\Gamma=\tau^2 I_p$ (hence ${\partial^{2}}/{\partial \Gamma^2}={\partial^{2}}/{\partial (\tau^2)^2}$).

\begin{corollary}\label{diagonalite.1}
Assume $\Gamma=\tau^2 I_p$ and $\theta=(\mu,\tau^2)\in\R^p\times \R_{+*}$. Then $\Psi_{\alpha}(\Gamma)=\tau^2 \Psi_{\alpha}$ with $\Psi_{\alpha}=\Sigma^{-\alpha/2} H  H^T \Sigma^{-\alpha/2}$
and
\begin{eqnarray*}
I_{{Y}_{\alpha}}(\mu) & = & H^T \Sigma^{-\alpha/2} \left(\alpha I_q + \tau^2\Psi_{\alpha}\right)^{-1} \Sigma^{-\alpha/2}  H, \\
 I_{{Y}_{\alpha}}(\Gamma) & = & \dfrac{1}{2}\sum\limits_{i=1}^q \left(\frac{\lambda^{\Psi_{\alpha}}}{\alpha + \tau^2\lambda^{\Psi_{\alpha}}}\right)^2.
\end{eqnarray*}
\end{corollary}

\subsubsection{Conditions for Fisher constraints}

From the previous results, some conditions for (\ref{eq:fisher-condition2}) can be derived. 
Note first that, from Theorem \ref{theo:fisher-info-gen},
\begin{eqnarray*}
 I_{g(X)}(\theta) - I_{Y^*}(\theta) & = & \left(
\begin{array}{ll}
I_{Y_0}(\mu) - I_{Y_1}(\mu) & {\bf 0} \\
{\bf 0} & I_{Y_0}(\Gamma) - I_{Y_1}(\Gamma)
\end{array}\right) \\
- \frac{1}{c}I_{g(X)}(\theta) + I_{Y^*}(\theta) & = & \left( \begin{array}{ll}
I_{Y_1}(\mu) - \frac{1}{c}I_{Y_0}(\mu) & {\bf 0} \\
{\bf 0} & I_{Y_1}(\Gamma) -  \frac{1}{c}I_{Y_0}(\Gamma)
\end{array}\right). 
\end{eqnarray*}
Then 
(\ref{eq:fisher-condition2}) is true if and only if the following conditions are simultaneously verified:
\begin{eqnarray}
I_{Y_0}(\mu) & \succeq & I_{Y_1}(\mu), \label{CN.1} \\
I_{Y_0}(\Gamma) & \succeq &  I_{Y_1}(\Gamma), \label{CN.3} \\
I_{Y_1}(\mu) & \succeq &  \frac{1}{c}I_{Y_0}(\mu), \label{CN.2} \\
I_{Y_1}(\Gamma) & \succeq &  \frac{1}{c}I_{Y_0}(\Gamma). \label{CN.4}
\end{eqnarray}
Hence necessary conditions (NC) for any of these inequalities are NC for (\ref{eq:fisher-condition2}). In the following results, we always assume that $H\Gamma H^T$ and $\Sigma$ are invertible. The result provided in Proposition \ref{going.further} highlights the difficulty of obtaining necessary and sufficient conditions (NCS) for (\ref{eq:fisher-condition2}) in general cases, as they would require NCS for differences of sums of products of orthonormal vectors that belong to different spaces, as well as prior knowledge on the ordering of pseudo-inverse matrices.  For this reason, several common simplified situations are examined in the following, in particular situations when the output $Y$ is unidimensional.

\begin{Prop}\label{NC.A}
Inequation (\ref{CN.1}) is always true.
\end{Prop}

\noindent Comparing $I_{Y_0}(\Gamma)$ and $I_{Y_1}(\Gamma)$ under the context of Proposition \ref{going.further} can be simplified if $\Psi_0(\Gamma)$ and $\Sigma^{1/2}$ are co-diagonalizable (namely if they share a common base of diagonalization). Especially, this is true if $\Sigma=\sigma^2 I_q$ or if $q=1$.

\begin{Prop}\label{NC.C}
Assume that $\Gamma$ is diagonal and $\Psi_0(\Gamma)=HH^T$ and $\Sigma^{1/2}$ are co-diagonalizable. Then Inequation (\ref{CN.3}) is true.
\end{Prop}

In other situations, with $q>1$ and $\Gamma$ still diagonal, denote $Y_{[i]}$ the $i$th-dimensional output of $Y\in\R^q$. Considering only this output dimension, then  $I_{{Y_{[i]}}_{0}}(\Gamma)  \succeq  I_{{Y_{[i]}}_{1}}(\Gamma)$.  
Note besides that (\ref{decomp.Igamma}) and (\ref{fisher.term.simp.2}) in Theorem \ref{theo:fisher-info-gen} and Proposition \ref{going.further} highlight the fact that, provided $\Sigma=\sigma^2 I_q$, then 
\begin{eqnarray*}
I_{Y_{\alpha}}(\Gamma) & = & \sum\limits_{i=1}^q I_{{Y_{[i]}}_{\alpha}}(\Gamma).
\end{eqnarray*}
Consequently, by transitivity of the Loewner order, the relation (\ref{CN.3}) still hold when $\Sigma=\sigma^2 I_q$.  

\begin{Prop}\label{NC.B}
Inequation (\ref{CN.2}) is true if and only if
\begin{eqnarray}
H \Gamma H^T & \succeq & \frac{1}{c-1} \Sigma.
\label{fisher.0.grosgris}
\end{eqnarray}
Furthermore,  a NC  for (\ref{CN.2}) is
\begin{eqnarray}
\|\Gamma\|_2  & \geq &  \frac{1}{c-1} \left\|(HH^T)^{-1} \Sigma\right\|_2.\label{cond0toto}
\end{eqnarray}
\end{Prop}

\noindent Note that this last result is actually a NC for (\ref{fisher.0.grosgris}), and it can be easily transformed into a sufficient condition (SC) for this latter inequation by replacing $\|\Gamma\|_2=\lambda^{\Gamma}_{\sup}$ by $\lambda^{\Gamma}_{\inf}$ (see the proof). \\

\begin{Prop}\label{NC.D}
Assume the conditions and notations of Corollary \ref{diagonalite.aa}. Denote $\rho =  \sqrt{\left\| \left(\Lambda_{1} \Lambda^T_{1}\right)^+  \left(\Lambda_{0} \Lambda^T_{0}\right) \right\|_2}$, $\gamma=\|\Sigma^{-1/2}\|^4_2\|\Lambda_1\|^{-2}_2 \|\Lambda_0\|^2_2$ and 
$U   =   \sqrt{c}\tilde{\Delta}_0 - \rho \tilde{\Delta}_1$. 
Then Inequation (\ref{CN.4}) is true if and only if
\begin{eqnarray}
{\rm{Tr}}(U\Gamma) & \geq & \rho 
\label{fisher.ncd}
\end{eqnarray}
and  if $c>\gamma^2$ a NC (\textit{resp.} SC) for (\ref{fisher.ncd}) is 
\begin{eqnarray*}
\|\Gamma\|_2  \ \ \ \ \text{(\textit{resp.} $ \tau^2_{\inf}$)} & \geq &  \frac{\rho}{\sqrt{c}\| \Lambda_0 \|^2_2 - \rho \| \Lambda_1 \|^2_2}.
\end{eqnarray*}
\end{Prop}

\paragraph{The case $\Gamma=\tau^2 I_p$.} Studying this simplistic situation allows to provide a clear rationale to select a value for $c$, allowing to connect the Sobol' and entropy-based constrained problems previously studied with the Fisher-based constrained problem. From Corollary \ref{diagonalite.1} one may derive the following results.

\begin{Prop}\label{all.res} 
Assume $\Gamma=\tau^2 I_p$ and $\theta=(\mu,\tau^2)\in\R^p\times \R_{+*}$. Denote $\Psi_{\alpha}=\Sigma^{-\alpha/2} H  H^T \Sigma^{-\alpha/2}$. Then a NC for 
\eqref{eq:fisher-condition2} is
\begin{eqnarray}
\tau^2 & \geq & \frac{1}{\sqrt{c}-1}\| \Psi_{1}\|^{-1}_2. 
\label{NC.SC.without.com}
\end{eqnarray}
\end{Prop}

\noindent Condition (\ref{NC.SC.without.com}) refines the general NC (\ref{cond0toto}), noticing  that $\| \Psi_1\|^{-1}_2=\left\|(HH^T)^{-1} \Sigma\right\|_2$ when $HH^T$ and $\Sigma^{1/2}$ are co-diagonalizable. Expression (\ref{NC.SC.without.com}) highlights a lower bound for the signal over noise ratio (MSNR ; see $\S$ \ref{fisher.carac}): necessarily,
\begin{eqnarray*}
\text{MSNR} & = & \tau^2 \| \Psi_{1}\|_2 \ \geq \ \frac{1}{\sqrt{c}-1}.
\end{eqnarray*}
It seems reasonable a priori to assume $\text{MSNR} \geq 1$ for ensuring that the information required for the inversion task is not overwhelmed by model noise. This leads to select $c$ such that
\begin{eqnarray}
c  \geq  4.
\end{eqnarray}
This lower bound for $c$ seems to make sense if we go back to Section \ref{linear-1D}.
Consider the example given by \eqref{eq:model_toy},
where $q=1$ and $H=a\in\mathbb{R}$; then Condition (\ref{NC.SC.without.com}) gives
$(\sqrt{c}-1)\tau^2 >{\sigma^2}/{a^2}$. 
For model \eqref{eq:model2}, where $q=1$ and $H=a\in\mathbb{R}^p$, this  gives
$(\sqrt{c}-1)\tau^2 >{\sigma^2}/{\parallel a \parallel^2}$. Choosing $c=4$, the conditions provided for  Sobol' or entropic senses are recovered: 
$\tau^2 \|a \|^2_2 \geq \sigma^2$. 
In the context of Proposition \ref{NC.D}, if $q=1$ and $c=4$, then $\rho=\sigma^2$, $\gamma=1$ and the NSC (\ref{fisher.ncd}) similarly becomes $\tau^2  \geq  \frac{\sigma^2}{\|a\|^2_2}$.


\subsection{Gaussian linearizable models}\label{variational.analysis}

\subsubsection{Theoretical approximation}

Let us go back to the more general case where $g$ is some deterministic function from $ \mathbb{R}^{p}$ to $ \mathbb{R}^q$, not necessarily linear. For nonlinear models, the intractability  of the likelihood function carries forward to the calculation of the Fisher information matrix, preventing us from establishing directly simple NC for Condition \eqref{eq:fisher-condition2}.
A commonly used approach to bypass this problem is to linearize the regression function. Consequently, $Y$ is still assumed to have its first and second order moments finite. \\

To avoid the technical difficulties raised by local linearity (see $\S$ \ref{locally.linear}), consider now that $g$ can be linearizable in a variational sense, at least on a given domain of interest.  In the context of this study, such an approximation -- replacing the output $Y$ by a tractable  approximation $\tilde{Y}$ -- should be considered relevant if it allows to carry the most part of information transferred by $g$. For this reason, the error resulting from the approximation should have a small impact on $\tilde{Y}$. Intuitively, this means that obtaining a Fisher constraint from the linearized model  should help the inversion of the nonlinear model,  which requires statistical estimation; this heuristic will be tested in practice in Section \ref{numeric}. \\

More formally, an information-theoretic rationale, usual in variational analysis, can be invoked to define $\tilde{Y}$ as the best approximation of $Y$ in a distributional sense using the Kullback-Leibler divergence to quantify the quality of the approximation \cite{Cover2012}. Recall that this divergence between two distributions $P$ and $Q$, with respective pdf $(f_p,f_q)$ that are absolutely continuous with respect to each other on a  common support $D$, is given by 
	$$\text{KL}\left(P\| Q \right)=\int_{D} f_p(x)\log \frac{f_p(x)}{f_q(x)}dx. $$

\begin{Prop}\label{prop4bis}
Denote $Y=g(X)$ and $Y^*=Y+\varepsilon$ with $\E[\varepsilon]=0$ and $X\sim{\cal{N}}_p(\mu,\Gamma)$ with $\Gamma\in S^{++}_p(\R)$. Denote ${\cal{L}}_{p,q}$ the space of homoscedastic Gaussian linear models from $\R^p$ to $\R^q$ and assume $q\leq p$. Define $U\sim{\cal{N}}(0,S)$ with $S$ independent of $X$ and for any $H\in\R^{q,p}$, define in ${\cal{L}}_{p,q}$ the linear model
\begin{eqnarray*}
\tilde{Y} & = & H(X-\mu) + \E[Y^*] + U.
\end{eqnarray*}
Assume $\mbox{Cov}(Y^*)\succeq \Sigma+S$ and denote $(f_Y,f_{\tilde{Y}_H})$ the pdf of $(Y,\tilde{Y})$, respectively. Then there exists  $H^*\in\R^{q,p}$ such that
\begin{eqnarray*}
f_{\tilde{Y}_H^*} & = & \arg\min\limits_{f \in {\cal{L}}_{p,q}} \text{KL}\left(f_Y,f\right)
\end{eqnarray*}
and $H^*$ is such that $H^*\Gamma H^{*T} + S  =  \mbox{Cov}(Y^*) - \Sigma$.
\end{Prop}

\noindent It follows from this last result, adapting Proposition \ref{NC.B}, that a NC for a meaningful linear approximation is
\begin{eqnarray*}
H^*\Gamma H^{*T} & = & \mbox{Cov}(Y^*) - \Sigma - S \ \succeq \ \frac{1}{c-1} \left(\Sigma + S\right)
\end{eqnarray*}
or similarly
\begin{eqnarray}
\mbox{Cov}(Y^*)  & \succeq & \frac{c}{c-1} \left(\Sigma + S\right). \label{observed.cov.condition}
\end{eqnarray}

Since $A\succeq B \Rightarrow {\rm{Tr}}(A)\succeq {\rm{Tr}}(B)$ for $(A,B)$ two square matrices, this implies directly the following result.

\begin{Prop}\label{final.NC}
A NC for a meaningful linear approximation of $g$ is
\begin{eqnarray}
{\rm{Tr}}\left(H^*\Gamma H^{*T}\right) & \geq & \frac{1}{c-1} {\rm{Tr}}\left(\Sigma + S\right). \label{final.nc.eq}
\end{eqnarray}
\end{Prop}

In the case where $\Gamma=\tau^2 I_p$, estimating $H^*$ appears useless to provide a condition on $\tau^2$, as stated by next Corollary.

\begin{corollary}\label{final.corollary}
Assume $\Gamma=\tau^2 I_p$. Then a NC for the linear approximation is (\ref{NC.SC.without.com}), replacing $\Psi_1$ by
\begin{eqnarray}
\Psi & = &  \left(\Sigma +S\right)^{-1/2} \left(\mbox{Cov}(Y^*) - S-\Sigma\right)\left(\Sigma +S\right)^{-1/2}.  \label{cote2bis}
\end{eqnarray}
\end{corollary}

\subsubsection{Computational aspects}

The previous results require to get consistent estimators of $\E[Y^*]$ and $\mbox{Cov}(Y^*)$. In this stochastic inversion context, since
	a sample of noisy observations $(y^*_i)_{1\leq i \leq n}$  is always available in practice, it can be used to compute the usual empirical estimators $(\bar{y}^*_n,{c}^*_n)$ 
\begin{eqnarray}
\bar{y}^*_n & = & \frac{1}{n} \sum\limits_{i=1}^n  {y}^{*}_i, \label{empiric.es1} \\
\bar{c}^*_n & = & \frac{1}{n} \sum\limits_{i=1}^n  \left({y}^{*}_i-\bar{y}^*_n\right)\left({y}^{*}_i-\bar{y}^*_n\right)^T.\label{empiric.es2}
\end{eqnarray}
While $H$ is not required for the computation of (\ref{cote2bis}), $S$ must be estimated too. Since this covariance matrix expresses the error between $g$ and the best linear model, a simple approach to estimating $S$  can be based on  using a numerical Design Of Experiments (DOE), or \emph{training sample}, $({x}'_{j},{y}^{'}_{j}=g({x}'_{j}))_{1\leq j \leq m}$ with $x'_j=(x'_{j,l})_{1\leq l \leq p}$ and $y'_j=(y'_{j,t})_{1\leq t \leq q}$, an usual tool for the preliminary exploration of computer models. \\

The choice of a DOE is part of an usual and general problem in numerical analysis of computer models, that depends on both the computational easiness of sampling with $g$ and the dimension of $X$. To avoid defining formal bounds on the numerical support of $X$, in a Bayesian informative setting where a prior $\theta\sim \pi(\theta)$ is assumed to be known, sampling uniformly the DOE within the highest predictive prior density region
$\{X; \ X  \sim  {\cal{F}}(\theta), \ \theta \sim\pi(\theta)\}$
is intuitive. But very generally, space-filling type approaches with projection regularity properties on the sampling subspaces, such as the maximin Latin Hypercube sampling, are preferred to avoid clusters and filling gaps. They require bounded domains that may be defined using extreme quantiles of the predictive prior density. See \cite{Damblin2013} for more details. 

The availability of the DOE allows to understand $S$ as the covariance of residuals in a multiple regression, consistently estimated by
\begin{eqnarray}
{s}^*_{n,m} & = & \left(W'_{m,n} - X'_m\hat{H}_{n,m}\right)\left(W'_{m,n}  - X'_m\hat{H}_{n,m}\right)^T \label{estimate.s}
\end{eqnarray}
where $W'_{m,n}\in\R^{m,q}$ and its $(i,j)$-element is $y'_{i,j}-\bar{y}^*_{n,j}$, $X'_m\in\R^{m,p}$ and  $\hat{H}_{n,m}\in\R^{p,q}$ solves the normal equation of the associated linear system, namely
\begin{eqnarray}
\hat{H}_{n,m} & = & \left(X^{'T}_m X'_m\right)^{-1} X^{'T}_m W'_{m,n}
\end{eqnarray}
provided the DOE allows $X^{'T}_m X'_m$ to be invertible. It has sense only if the resulting approximation of the left term of  (\ref{cote2bis}) is positive definite, hence if
\begin{eqnarray}
 c^*_n -  {s}^*_{n,m} & \succ & \Sigma. \label{cond.covar}
\end{eqnarray}

\section{Numerical experiments}\label{numeric}

\subsection{General framework}\label{gen.ex}

This section illustrates the effect of constraining the variance of $X$ in an usual Bayesian setting. Consider the  stochastic inverse problem (\ref{inversion.problem}) and the general prior form advocated in \cite{Fu2015} and \cite{Fu2016}:
\begin{eqnarray}
\mu|\Gamma & \sim & {\cal{N}}_p\left(\mu_0,\Gamma/\alpha\right), \label{prior1} \\
\Gamma & \sim &  \mathcal{IW}_p(\Upsilon,\nu).\label{prior2}
\end{eqnarray}
where $\mathcal{IW}$ denotes the inverse Wishart distribution on $S^{++}_p(\R)$, \textcolor{black}{and where the prior hyperparameters $(\mu_0,\alpha,\Upsilon,\nu)$ are fixed.} 
This choice is often made for conjugacy reasons. Indeed, given an observed sample $\mathbf{y}^*=(y^*_1,\ldots,y^*_n)$, a Gibbs sampler adapted from  \cite{Fu2015} is a natural approach to generate samples from the posterior distribution $\pi(\mu,\Gamma|y^*_1,\ldots,y^*_n)$ \cite{TL95}, that implements a data augmentation scheme and benefit from conditional form stability.   Given initializing values $(\mu^{[0]}, \Gamma^{[0]})$ for the parameters, and initialized values  $\mathbf{X}^{[0]}$ for the missing variables $\mathbf{X}=(X_1,\ldots,X_n)$, generate the following Markov chain  for $r=0,1,2,\dots$: 

\setcounter{equation}{3}
\noindent\rule{\linewidth}{.5pt}
\begin{eqnarray}
\text{1.} \ \ \Gamma^{[r+1]} & \sim & \mathcal{IW}_p\left(\Upsilon+\sum_{i=1}^n (\mu^{[r]}-X_i^{[r]})(\mu^{[r]}-X_i^{[r]})^T+\alpha(\mu^{[r]}- \mu_0)\times\right. \nonumber \\
& & \ \hspace{1cm} \left. (\mu^{[r]}-\mu_0)^T, \,\nu+n+1\right), \label{sim.gamma} \\  
\text{2.} \ \ \mu^{[r+1]} & \sim & \mathcal{N}_p\Big(\frac{\alpha}{n+\alpha}\mu_0+\frac{n}{n+\alpha}\overline{\mathbf{X}_n^{[r]}},\, \frac{\Gamma^{[r+1]}}{n+\alpha}\Big) \ \ \text{
 with
$\overline{\mathbf{X}_n^{[r]}}=\frac{1}{n} \sum\limits_{i=1}^n X_i^{[r]} $}, \nonumber \\ 
\text{3.} \ \ \mathbf{X}^{[r+1]}  & \sim & \tilde{\pi}\left(\mathbf{X}^{[r+1]}|g,\mathbf{y}^*,\Sigma,\mu^{[r+1]}, \Gamma^{[r+1]} \right), \nonumber
\end{eqnarray}
this latter conditional density being proportional to
\setcounter{equation}{4}
\begin{eqnarray}
 \exp\Big\{-\frac{1}{2} \sum_{i=1}^n \Big[ \left(X^{[r+1]}_i-\mu^{[r+1]}\right)^T\left(\Gamma^{[r+1]}\right)^{-1}\left(X^{[r+1]}_i-\mu^{[r+1]}\right) \nonumber \\ 
  + \left(y^*_i-g(X^{[r+1]}_i)\right)^T \Sigma^{-1}\left(y^*_i-g(X^{[r+1]}_i)\right)\Big]\Big\}. \label{density.of.x}
\end{eqnarray}
\rule{\linewidth}{.5pt}

Since (\ref{density.of.x}) does not belong to a closed form family of density functions, a Metropolis-Hasting (MH) step can be used to  simulate $\mathbf{X}^{[r+1]}$ from its full conditional distribution \cite{TL95}. Following the advices of \cite{Fu2015} resulting from many numerical tests, a multivariate random walk sampling candidates $\mathbf{X}_n^{[r]}$ can be an efficient instrumental distribution for this MH step.  
An additional gain of choosing the prior (\ref{prior1}-\ref{prior2}) is that the meaning of hyperparameters $(\mu_0,\alpha,\Upsilon,\nu)$ is straightforward, as described in \cite{Fu2016}. 

\subsection{Motivating case-study}

In this later article, the authors considered a hydraulic engineering model $g(X,d)$ previously used in \cite{bastos09} where $d$ is a known water discharge value of a river, $X$ is a unknown, two-dimensional vector of roughness coefficients for the riverbed and floodplain, respectively, and $Y$ is the corresponding bivariate vector of water heights (levels) in each part of the geographical area. A simplified (while rather good) approximation of $g$ is the following: 
\begin{eqnarray}
Y = \left(X_2 + \left(\dfrac{d \cdot \sqrt{5000}}{300 \cdot \sqrt{55 - X_2}} \cdot X_1\right)^{0.6},  \dfrac{d^{\ 0.4} \cdot X_1^{0.6} \cdot (55 - X_2)^{0.3}}{300^{0.4} \cdot 5000^{0.3}}\right). \label{eq:h-model}
\end{eqnarray}

In the remainder of this section, numerical experiments are conducted to compare the matching between the posterior predictive distribution of the missing values $X$ computed using the previous algorithm with a chosen simulation distribution, in settings where the prior  (\ref{prior2}) is constrained or not by the condition $\Gamma\in\Omega^{(k)}_j$ ($j\in\{A,B,C\}$). ``Constrained by" means that Step 1 of the algorithm  is replaced by a MH step where:
\begin{itemize}
\item the instrumental distribution for $\Gamma$ is proportional to
\begin{eqnarray*}
& \mathcal{IW}_p\left(\Upsilon+\sum_{i=1}^n (\mu^{[r]}-X_i^{[r]})(\mu^{[r]}-X_i^{[r]})^T+\alpha(\mu^{[r]}- \mu_0)\times \right. \nonumber \\
&  \ \hspace{3cm} \left. (\mu^{[r]}-\mu_0)^T, \,\nu+n+1\right) \1_{\left\{\Gamma \in \Omega^{(k)}_j\right\}} , 
\end{eqnarray*}
which can be simulated by simple acceptation-rejection;
\item the Metropolis ratio does not involve the constraint $\1_{\left\{\Gamma \in \Omega^{(k)}_j\right\}}$. 
 \end{itemize}
 This appears necessary since, with $\Gamma\in\Omega^{(k)}_j$, the conditional posterior of $\Gamma$ given the other parameters is no longer an Inverse Wishart distribution. It is likely that better sampling techniques can be used, but they were not investigated in this article.

The simulation data $X$ are generated according to $f_s(X)\equiv{\cal{N}}(\mu_s,\Gamma_s)$ where $\mu_s=(30,50)$ and $\Gamma_s={diag}(5^2, 1)$, while several values of $d$ are considered, sampled from $f_d \equiv\text{Gumbel}(1013, 100)$, that define several instances of $g$. Such choices lead to physically plausible observable values 
\setcounter{equation}{5}
\begin{eqnarray}
Y^{(k)*}_i & = & g\left(X_i,d^{(k)}_i\right) + \varepsilon_i, \ \ \ \text{$i\in\{1,\ldots,n\}$, $k\in\{1,\ldots,M\}$,} \label{simulation.xy}
\end{eqnarray}
in a 
 real context where the model is used to invert roughness coefficients of a French river \cite{Fu2016}. Following the magnitudes provided by these authors, the covariance matrix of noise $\varepsilon$ is chosen as $\Sigma = {diag}(1^2, 0.1^2)$ 
  and the chosen prior hyperparameters for  (\ref{prior1}-\ref{prior2}) are summarized on Table \ref{table:prior1}. 
 
 \begin{table}[hbtp]
            \caption{Prior hyperparameters.}
            \label{table:prior1}
    \begin{tabular}{cc}
        \hline
        $\mu_0$ & (35, 49) \\
        $\alpha$ & 1 \\
        $\nu$ & 5 \\
        $\Upsilon$ & $2\begin{pmatrix}
        7.5^2 & 0 \\
        0 & 1.5^2
        \end{pmatrix}$ \\
        \hline
    \end{tabular}
\end{table}

\subsection{Models and prior constraints}\label{linear.surrg}

In this framework where $g$ is not linear, a class of $M$ simplified linear models $\hat{g}^{(1)},\ldots, \hat{g}^{(M)}$, or \emph{surrogates} approximating $g$ was produced as follows, for the experimental needs. Following approaches advised by \cite{Iooss.Boussouf.ea2010} in a broader framework of Gaussian meta-modelling,  a maximin Latin Hypercube design of $N=2700$ values of $X$ was sampled within the cubic domain ${\cal{A}}=[20,50]\times[30,60]$ instead of a quasi-uniform one, and the corresponding values $Y^{(k)}=g(X,d^{(k)})$ were computed, given a sampled value of $d^{(k)}\sim f_d$, for $k\in\{1,\ldots,M\}$ with $M=200$. $M$ multivariate linear (surrogate) models were fitted from the $M\times N$ couples $(x^{(k)}_i,y^{*(k)}_i)_{i,k}$, such that 
\setcounter{equation}{7}
\begin{eqnarray}
y^{*(k)}_i & = & \hat{g}_M(x^{(k)}_i) + \varepsilon_i, \nonumber \\
 \hat{g}_M(x^{(k)}_i) & = &   y^{(k)}_i + u^{(k)}_i \ = \  a^{(k)} x^{(k)}_i + b^{(k)} + u^{(k)}_i  \label{linear.class}
\end{eqnarray}  
where each  $u^{(k)}_i$ was estimated as the bivariate centered Gaussian noise ${\cal{N}}(0,S^{(k)})$, approximating the model error between $g$ and $\hat{g}^{(k)}$,  where each $S^{(k)}$ is estimated from (\ref{estimate.s}), with
\begin{eqnarray*}
\bar{y}^{(k)*}_n  \simeq   \left(\begin{array}{l} 
52.13  \\
1.76 
\end{array}\right) & \text{and} & \bar{c}^{(k)*}_n   \simeq  \left(\begin{array}{ll} 
2.6 & -0.14 \\
-0.14 & 0.038
\end{array}\right)
.     
\end{eqnarray*}
Details on computations and diagnostics are provided in Appendix \ref{linear.append}. On average, the sum of $\Sigma$ and $S^{(k)}$ was found to be 
\begin{eqnarray*}
 \E\left[\Sigma+S^{(k)}\right] &=&
\begin{pmatrix} 
    1.008 & -0.003 \\
    -0.003 & 0.012
    \end{pmatrix}
\end{eqnarray*}
while the average values of matrix $a^{(k)}$ and  vector of intercepts $b^{(k)}$ are, respectively
\begin{eqnarray*}
\E\left[a^{(k)}\right] &=& \begin{pmatrix} 
    0.03 &  0.029 \\
     1.10 & -0.064 
    \end{pmatrix} \ \ \ \ \ \text{and} \ \ \ \E\left[b^{(k)}\right] \ = \ \left(\begin{array}{l} -2.53 \\  3.88 \end{array}\right). 
\end{eqnarray*}
For the $M$ models, it was checked that $\left|\Sigma\right| = 0.01 \leq \left|\Sigma+S^{(k)}\right| \sim 0.012 \leq  \left| a^{(k)} \Gamma_s a^{(k) T}\right| \sim 0.027$, ie. the uncertainty in $X$ is of primary influence in the uncertainty in $Y$, according to the entropy-based notion of uncertainty expressed by Definition (\ref{entropy}). Roughly speaking, each linear surrogate transmits a signal that is not overwhelmed by the noise. Besides, empirical condition (\ref{cond.covar}) was respected for all simulated samples, with 
\begin{eqnarray*}
\left| \bar{c}^{(k)*}_n - \Sigma - S^{(k)}\right| & > & 0.01. 
\end{eqnarray*}  
Finally, two situations of stochastic inversion are studied:
\begin{description}
\item[Situation 1.] We focus only on inverting the roughness coefficient from the water level on the riverbed, namely the first dimension of $Y$ (named $Y_{[1]}$), as collecting regular observations $Y^*_{[1]}$ is the most traditional approach to feed such hydraulic problems. Then relevant constraints for $\Gamma$ arise from considering $q=1$. 
\item[Situation 2.] We focus on solving the inverse problem considering the whole set of two-dimensional observations.  
\end{description}
Then, several constraints on $\Gamma$, associated to $g$ or to each $\hat{g}^{(k)}$ are now defined for the experiments. They  simply take the form of indicator constraints for defining the supporting subset of $S^{++}_p(\R)$ for $\Gamma$ in (\ref{prior2}). For all $k\in\{1,\ldots,M\}$, consider: 
\begin{eqnarray*}
\Omega^{(k)}_A\left(\Sigma\right) & = & \left\{\Gamma\in S^{++}_p(\R), \ \left| a^{(k)} \Gamma a^{(k) T}\right| \geq \left|\Sigma\right|   \right\}, \\
 \Omega^{(k)}_B & = & \Omega^{(k)}_A\left(\Sigma+S^{(k)}\right), \\ 
  \Omega^{(k)}_C & = & \left\{\Gamma\in S^{++}_p(\R), \   \| \Gamma \|_2 \geq \frac{1}{3} \left\| \left(a^{(k)}a^{(k) T}\right)^{-1} \Sigma\right\|_2  \right\}, \\
    \Omega^{(k)}_D & = & \left\{\Gamma\in S^{++}_p(\R), \   \rm{Tr}\left(\left(2\tilde{\Delta}^{(k)}_0 - \rho^{(k)} \tilde{\Delta}^{(k)}_1\right)\Gamma\right) \geq   \rho^{(k)}\right\}, \\
    \Omega^{(k)}_E & = & \left\{\Gamma\in S^{++}_p(\R), \   {\rm{Tr}}\left(a^{(k)}a^{(k) T}\right) \geq    \frac{1}{3} {\rm{Tr}}\left(\Sigma +S^{(k)}\right)\right\}.     
\end{eqnarray*}
Domain $\Omega^{(k)}_A$ appears as a soft  entropic prior constraint on $\Gamma$ coming from (\ref{eq:model2ter1}), that assumes negligible model error with respect to a linear surrogate. It can be rightfully tested over models $\hat{g}^{(k)}$ and $g$. Domain $\Omega^{(k)}_B$ is a more stringent condition that accounts for this model error, that should be tested over model $g$ only. Domain $ \Omega^{(k)}_C$ comes from (\ref{fisher.0.grosgris}), using $c=4$. It is appropriate for linear surrogates  $\hat{g}^{(k)}$ but it can be tested too over model $g$. Domain $ \Omega^{(k)}_D$, coming from (\ref{fisher.ncd}) with $c=4$, can be tested in a similar way but only makes sense in Situation 1 (univariate output of interest), replacing $\Sigma$ by its first diagonal component $\Sigma_{[1,1]}$. The composite constraint $\Omega^{(k)}_C \cap  \Omega^{(k)}_D$ can be also be rightfully tested in Situation 1. Finally, $ \Omega^{(k)}_E$ comes from (\ref{final.nc.eq}) with $c=4$ and is applicable only on model $g$. Tested models with their prior domains (or indicator constraints) for  $\Gamma$  are summarized on Table \ref{summary}. \\

\begin{table}[hbtp]
\caption{Models and associated prior covariance constraints tested.}
\label{summary}
\begin{tabular}{llll}
\hline
Indicator & Models & Prior domain  & Output dimensionality $q$ used for inversion \\
\hline
& {\it Linear surrogates} & \\
\cline{2-2}
& \\
$L_1$ & $\hat{g}^{(k)}(X)$ & $S^{++}_p(\R)$ (no prior constraint) & $q\in\{1,2\}$ \\[2mm]
$L_2$ & $\hat{g}^{(k)}(X)$ & $\Omega^{(k)}_A(\Sigma)$ & $q\in\{1,2\}$   \\[2mm]
$L_3$ & $\hat{g}^{(k)}(X)$ & $\Omega^{(k)}_C$ & $q\in\{1,2\}$  \\[2mm]
$L_4$ & $\hat{g}^{(k)}(X)$ & $\Omega^{(k)}_D$ & $q=1$  \\[2mm]
$L_5$ & $\hat{g}^{(k)}(X)$ & $\Omega^{(k)}_C \cap  \Omega^{(k)}_D$ & $q=1$  \\[2mm]
& \\
& {\it Nonlinear models} & \\
\cline{2-2}
& \\
$T_1$ & $g\left(X,d^{(k)}\right)$ & $S^{++}_p(\R)$ (no prior constraint) & $q\in\{1,2\}$  \\[2mm]
$T_2$ & $g\left(X,d^{(k)}\right)$ & $\Omega^{(k)}_A(\Sigma)$ & $q\in\{1,2\}$  \\[2mm]
$T_3$ & $g\left(X,d^{(k)}\right)$ & $\Omega^{(k)}_B$ & $q\in\{1,2\}$  \\[2mm]
$T_4$ & $g\left(X,d^{(k)}\right)$ & $\Omega^{(k)}_C$ & $q\in\{1,2\}$  \\[2mm]
$T_5$ & $g\left(X,d^{(k)}\right)$ & $\Omega^{(k)}_D$ & $q=1$  \\[2mm]
$T_6$ & $g\left(X,d^{(k)}\right)$ & $\Omega^{(k)}_C \cap  \Omega^{(k)}_D$ & $q=1$  \\[2mm]
$T_7$ & $g\left(X,d^{(k)}\right)$ & $\Omega^{(k)}_E $ & $q\in\{1,2\}$  \\[2mm]
\hline
\end{tabular}
\end{table}

\subsection{Experiments}

For $k\in\{1,\ldots,M=200\}$, a $n-$sample $Y^{(k)}$  of vectors is generated from (\ref{simulation.xy}). For each sample and each couple of model and prior constraints, an instance of the Metropolis-Hastings-within-Gibbs (or \emph{hybrid} Gibbs) algorithm described in $\S$ \ref{gen.ex} is run, adding each prior indicator constraint for $\Gamma$ to the simulation step (\ref{sim.gamma}).  Three Markov chains are run in parallel for each scenario, and the convergence towards each stationary, marginal posterior distribution of $X$ is assessed by  satisfying the usual condition on Gelman's statistic $\hat{R}_G < 1.02$ \cite{Gelman1998} (more details are provided on Appendix \ref{gelman}). Convergence is usually observed from 6000 iterations and chains are run such that, after a decorrelation step, posterior samples contain 5000 vectors.

\begin{figure}[hbtp]
\centering
\includegraphics[scale=0.45]{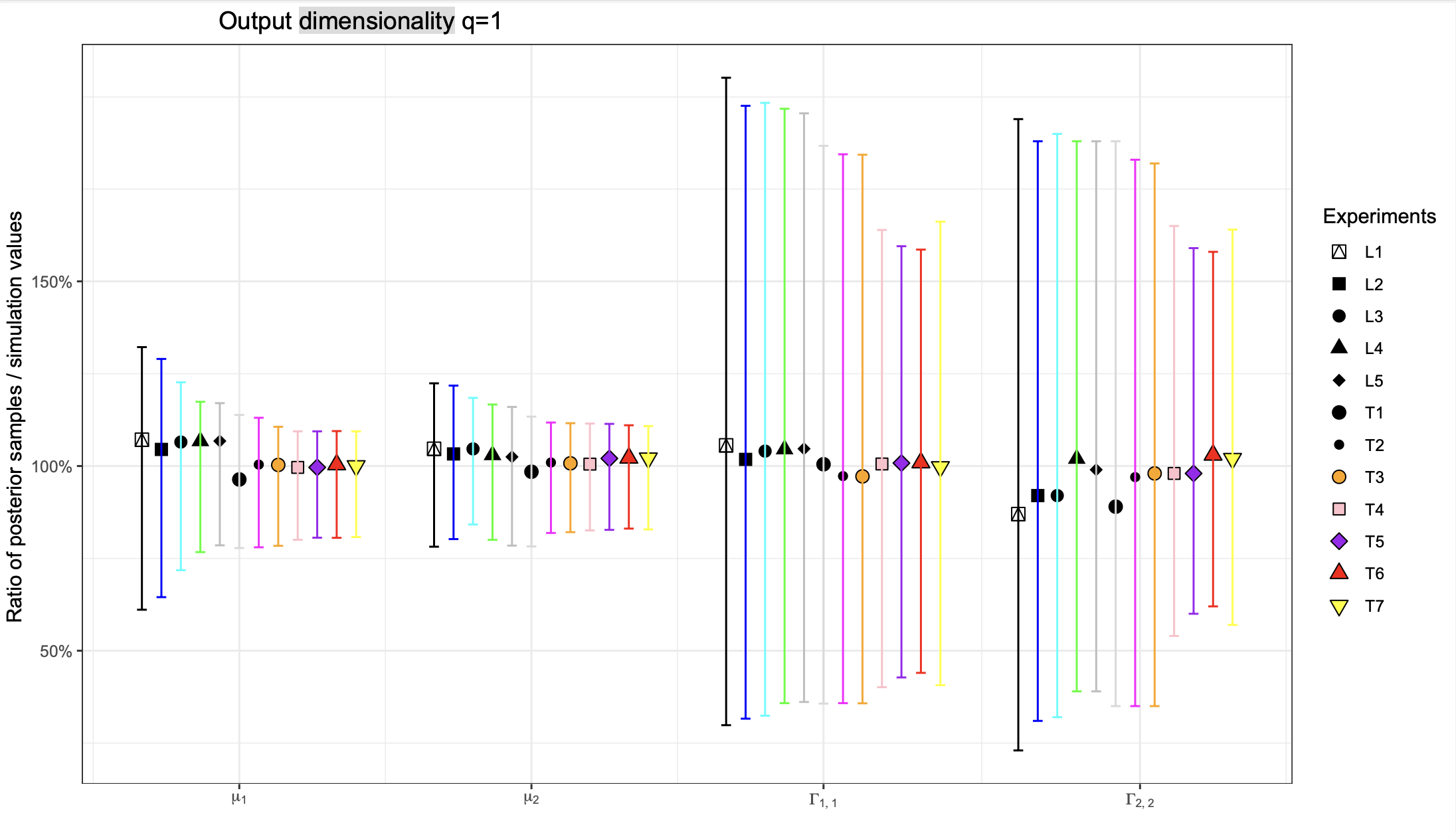}
\caption{Summary of posterior relative means and 90\% credibility intervals with respect to simulation values, on average on the values of $d$, for $n=30$. The inversion algorithm considers {\bf Situation 1} ($q=1$). Means are perfect estimates of simulation values when the $y-$axis values reach 100\%. }
\label{fig:posterior-1}  
\includegraphics[scale=0.46]{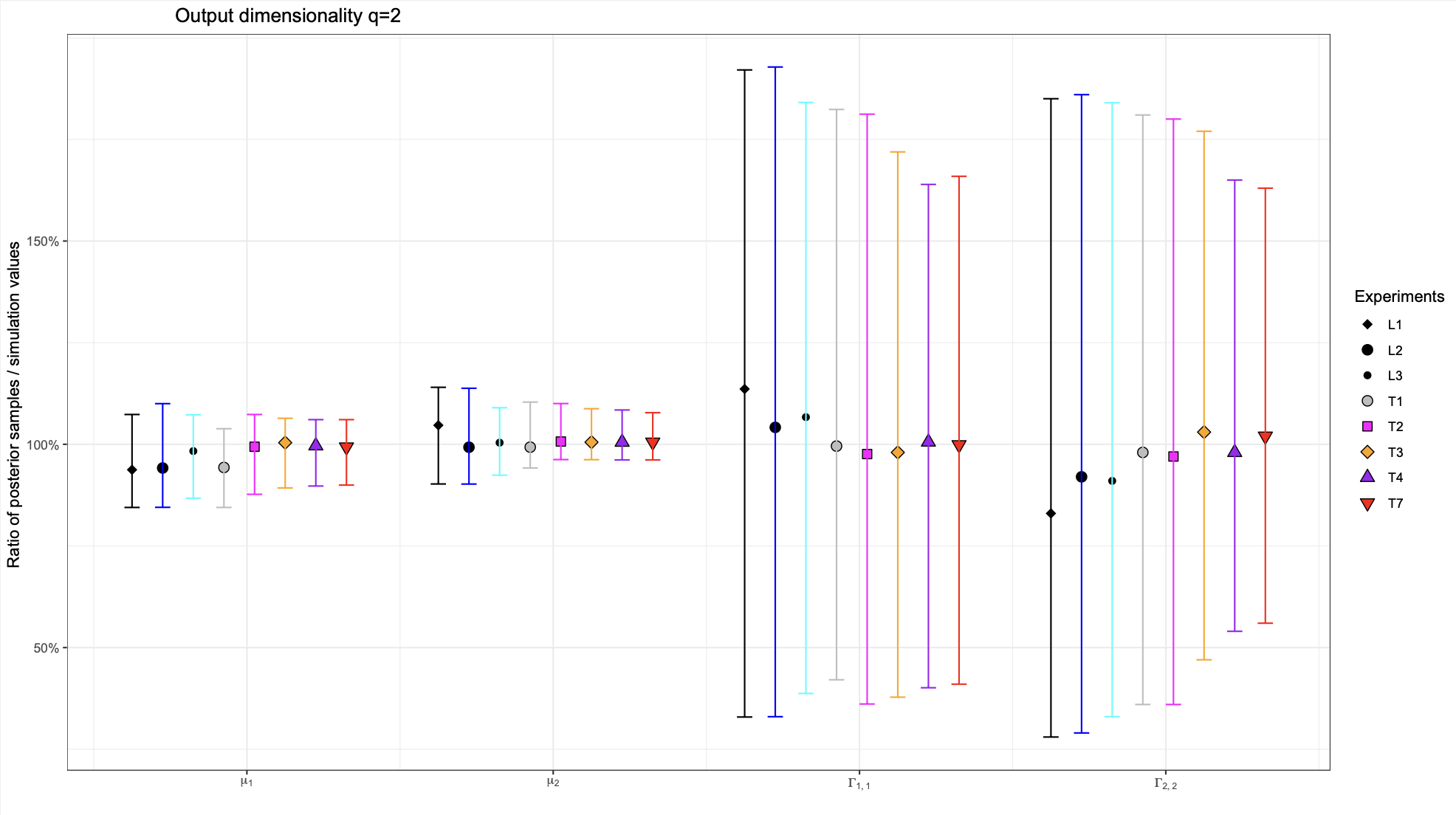}
\caption{Summary of posterior relative means and 90\% credibility intervals with respect to simulation values, on average on the values of $d$, for $n=30$. The inversion algorithm considers {\bf Situation 2} ($q=2$). Means are perfect estimates of simulation values when the $y-$axis values reach 100\%. }
\label{fig:posterior-2}  
\end{figure}

By design in \cite{Fu2016}, the prior remains weakly informative and in agreement with the generated data (in the sense of \cite{Bousquet2008}). 
The ability of the algorithm to estimate correctly the parameters, as well as the effect of constraints to shrink the highest posterior regions around the true values, can be appreciated on Figure \ref{fig:posterior-1} for Situation 1 ($q=1$) and Figure \ref{fig:posterior-2} for Situation 2 ($q=2$). While these figures show that these effects may be non-negligible (both on $\mu$ and $\Gamma$ through their entanglement in the estimation algorithm), their impact on the relevance of the posterior predictive distribution of $X$ can be more finely addressed as follows. Each marginal posterior on $X$, denoted by its density function $f_{S_j}(x|\hat{y}^{*(k)},d^{(k)})$ where $S\in\{L,T\}$ and $i\in\{1,\ldots,6\}$, provides an approximation of the simulation distribution $f_S$, and it can be compared to other posteriors through error indicators $e_{S_j}$ defined with respect to this latter distribution.  The  Kullback-Leibler divergence between the simulation distribution $f_s(X)$ and each marginal posterior is used here for defining this indicator:
\begin{eqnarray*}
e^{(k)}_{S_j} & = & KL\left(f_s(.) \| f_{S_j}(.|\hat{y}^{*(k)},d) \right), \\
& = & -H(f_s) -  \E_{f_s}\left[\log f_{S_j}(.|\hat{y}^{*(k)},d) \right]
\end{eqnarray*}
where $H(f_s)$ is the entropy of $f_s$, provided by (\ref{eq:entropy_gaussian}). The second term is estimated by Monte Carlo: for a large value $Q=100,000$, 
\begin{eqnarray*}
\E_{f_s}\left[\log f_{S_j}(.|\hat{y}^{*(k)},d) \right] & \simeq & \frac{1}{Q}\sum\limits_{l=1}^Q \log \hat{f}_{S_j}(x_j|\hat{y}^{*(k)},d_k)
\end{eqnarray*}
 where  $x_j\overset{iid}{\sim} f_s$ and $\hat{f}_{S_j}$ is a two-dimensional Gaussian-kernel density estimator of  ${f}_{S_j}$ defined from posterior sampling, defined on the bounded domain ${\cal{A}}$. \\

\begin{figure}[hbtp]
\centering
\includegraphics[scale=0.38]{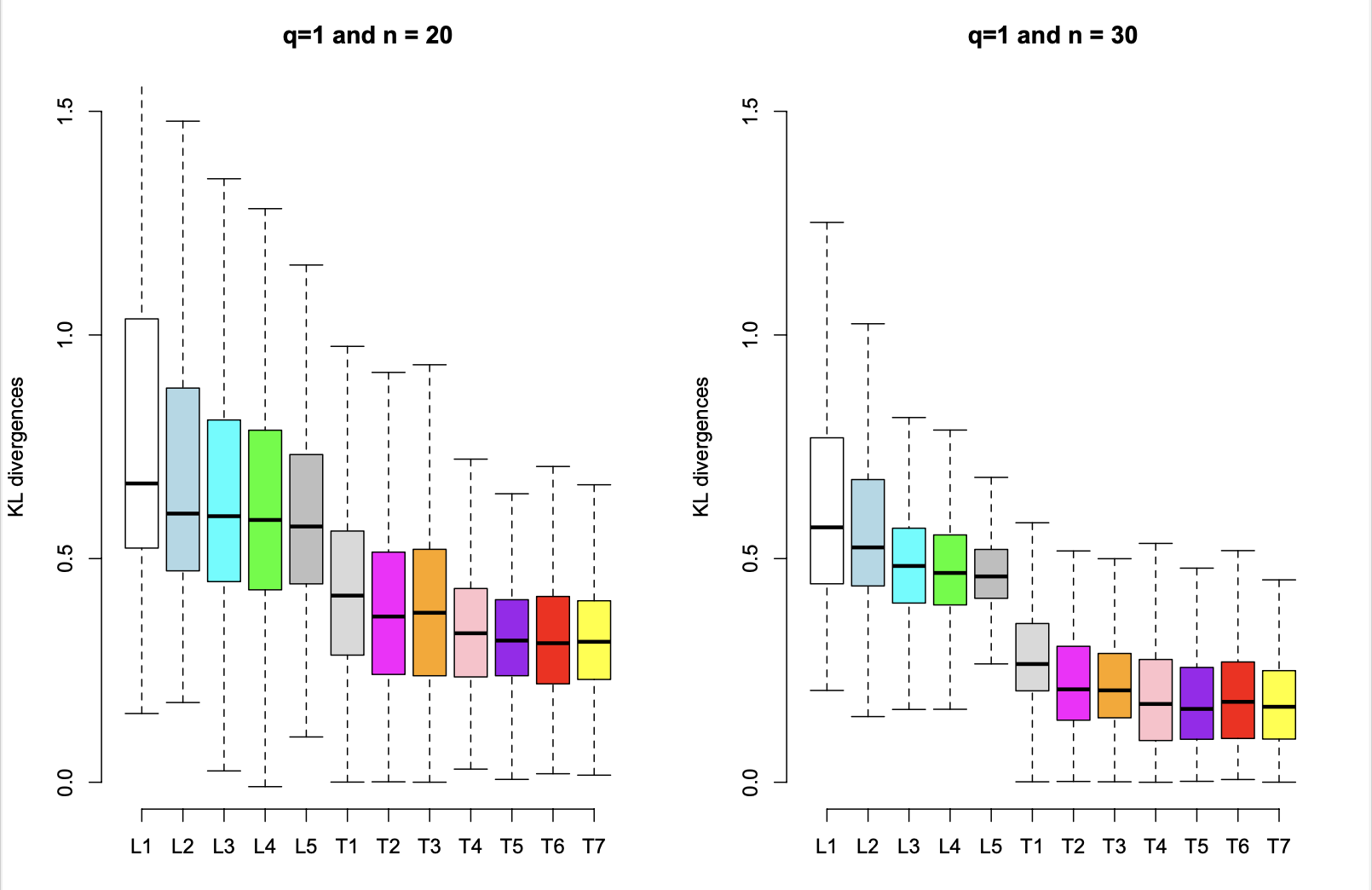}
\caption{Boxplots of errors $e^{(k)}_{S_j}$ (KL divergences) between classes of models and the simulation model (\textbf{Situation 1}: $q=1$), for two different data sizes.}
\label{boxplots1}
\includegraphics[scale=0.38]{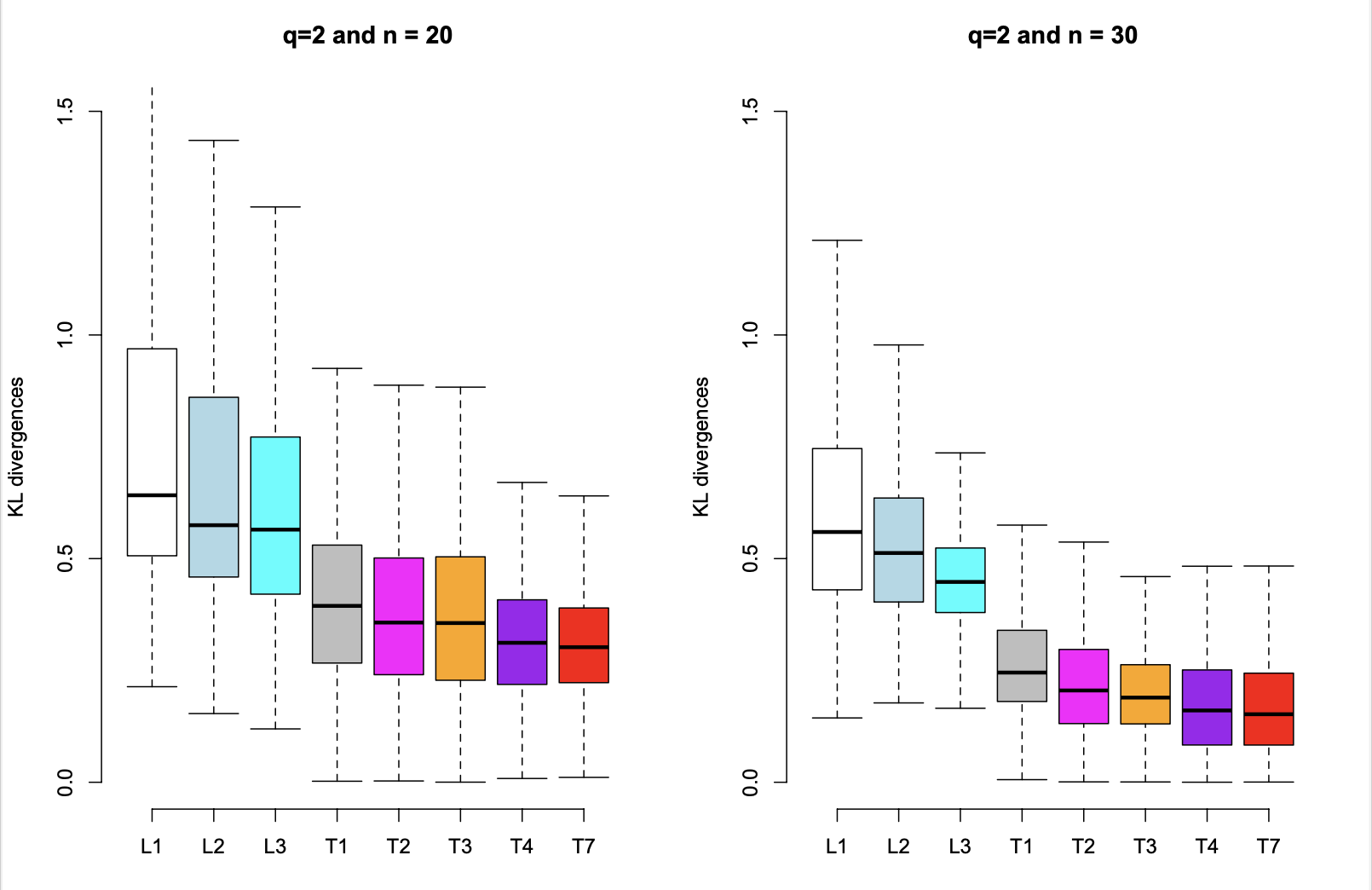}
\caption{Boxplots of errors $e^{(k)}_{S_j}$ (KL divergences) between classes of models and the simulation model (\textbf{Situation 2}: $q=2$), for two different data sizes.}
\label{boxplots2}
\end{figure}

Boxplots summarizing  how the errors $e^{(k)}_{S_j}$ are distributed are displayed on Figures \ref{boxplots1} and \ref{boxplots2} . As it could be expected, the posterior distributions estimated from the linear surrogates (L1 to L5)  have a larger error relative to the target distribution than the posterior distributions produced from the sampling model. When the dimension of output observations increases (from Situation 1 to Situation 2), posterior predictive distributions appear slightly closer to the true simulation distribution, as it could be expected again. Since all linear surrogates $\hat{g}^{(k)}$ are well defined according to (\ref{eq:model2ter1}), and because it accounts for model error in the prior constraint, it was expected that the posterior $L_3$  be closer to the simulation distribution ${\cal{N}}(\mu_s,\Gamma_s)$ than $L_1$ and   $L_2$. This behavior can scarcely be observed on these plots, because of the weakness of this error. While the proximity of any  posterior based on linearity approximation with the target distribution diminishes with the size $n$, it can nevertheless be observed that, on this example, the constraints built from linear approximations improve the separation of signal and noise within the stochastic inversion algorithm when it involves the true sampling model. The constraint resulting from a variational approximation seems to be as effective as the other constraints, while the latter are not completely relevant for all situations.  More generally, on this example inserting prior constraints help on average getting better recovery of the predictive posterior distribution and producing posterior estimates that fit better with the features of the underlying true distribution of $X$.

\section{Discussion}

Using tools of sensitivity analysis to formalize ideas about the relevance of solutions of stochastic inverse problems, that incorporate model noise $\varepsilon$, seems fruitful. Assuming a Gaussian input $X\sim{\cal{N}}_p(\mu,\Gamma)$ and considering linear and linearizable problems, this formalization leads to inequalities involving features of $\Gamma$ and $\varepsilon$, while offering (to our knowledge) a novel view on the degradation of Fisher information in noisy situations. These inequalities can be interpreted as prior constraints on $\Gamma$. Initial experiments show that, when incorporated into inversion algorithms, these constraints can improve Bayesian problem-solving (estimating the distribution of $X$) by reducing the size of the search space. 

However, this study remains a preliminary analysis that needs to be further developed, while paving the way for multiple improvements. If we first consider the practical use of such constraints, we notice that most inequalities on features of $\Gamma$ derived in this paper appear to be lower bounds for these features. This means that, when these inequalities are imposed as ``hard" prior constraints in the inversion process (in the sense these constraints are introduced using indicator functions), there is a risk that $\Gamma$ be conservatively overestimated. Although additional repeated experiments have not highlighted this phenomenon, it remains possible that the Bayesian estimation of $\Gamma$ will become biased when the number $n$ of available data grows to infinity. Such a behavior is not unusual in Bayesian settings when the prior is too strong. While this potential behavior  might be considered as beneficial for simulation studies related to risk and reliability analysis \cite{aven2016use}, this  limitation prompts us to define smoother constraints. For example, a hard constraint introduced using  Heaviside step functions as  $t(\Gamma)=\1_{\|\Gamma|\geq D_{g,\varepsilon}}$ could be smoothed using a bump transition function \cite{fry2002smooth}:
\begin{eqnarray*}
\hat{t}_m(\Gamma)  =  \left\{\begin{array}{l}
0 \ \ \text{if} \ \   |\Gamma|\leq D_{g,\varepsilon}(1-1/m) \\
1 \ \ \text{if} \ \     |\Gamma|\geq D_{g,\varepsilon}(1+1/m) \\
\frac{1}{2}\left(1 + \tanh\left[m^2\frac{\|\Gamma\|/D_{g,\varepsilon} -1}{1-m^2\left(\|\Gamma\|/D_{g,\varepsilon} -1\right)^2}\right]\right) \ \ \ \text{else.}
\end{array} \right. & \xrightarrow[]{m\to \infty} & t(\Gamma)
\end{eqnarray*}
An in-depth study of this approach would certainly be beneficial for improving the inferential framework, and also for integrating other types of prior information about $\Gamma$ in the same way. 
It is likely that the use of other distributions involving covariance structures would lead to similar constraints on this covariance, particularly due to the moment-matching properties resulting from KL minimization in the variational framework (e.g., exponential, hyperspherical and Mise-Fisher families ; see \cite{Kurz2016} and references therein). This framework could be studied more extensively, for example by replacing KL by a more general class of divergences, as the $\alpha-$divergences proposed by \cite{minka2005divergence} that could lead to tractable results. Furthermore, improvements to these constraints (e.g., obtaining necessary and sufficient conditions on meaningfulness by handling Fisher information) could potentially be achieved by leveraging recent advances in trace inequalities and connections between matrix partial ordering and statistical applications \cite{Golubic2020}.   

Concurrently, a research axis could involve adapting the most modern computational methods (e.g, \cite{Nagel2016,Nemeth2019}), which can accelerate the standard Monte Carlo approaches used in this paper:  ideally, instrumental distributions on $\Gamma$ should be supported on the intersection of the cone of semidefinite positive matrices with the subspaces defined by these constraints.

If we then seek to relax the assumptions on how to formalize the influence, as the choice of Fisher information is limited to models for which it is at least not singular,  the most modern concepts in GSA could be useful, such as those based on Bregman divergences \cite{DaVeiga13} or cooperative game theory. In particular, Shapley values \cite{Owen2014} or Proportional Marginal Effects  \cite{herin2024proportional} are variance-based coalition tools that are useful for avoiding the hypothesis of independence between $X$ and  $\varepsilon$. Focusing on entire distributions as uncertainty measures, Hilbert-Schmidt sensitivity indices \cite{DaVeiga13} and  probabilistic sensitivity measures recently defined by \cite{Antoniano2019} offer even greater generality, while being estimable in cases where the number of simulation data is limited --- an additional benefit for inverse problems that often incur high computational costs. In parallel, other approaches to  stochastic inversion \cite{Soize2013} could also be studied from the point of view proposed in this article, and adapted. In particular, spectral approaches such as polynomial chaos expansion \cite{Lemaitre2010,Yang2017} offer efficient calculation methods for hierarchizing uncertainties.

Beyond these technical improvements,  this paper highlights the existence of deep connections between the common UQ concepts of influence and concepts related to Hadamard's well-posed problems, such as $L^2$ conditioning, which are themselves tied to regularization concepts. Well-posed problems are generally studied in the context of deterministic inversion, where the values of $X$ are to be found. It is possible that this link appears solely due to the linearity assumptions used, but this remains to be explored in more general situations where $g$ is nonlinear. The results obtained by applying variational analysis suggest that exploring so-called meaningful inversion in nonlinear cases will rely on exploiting a combination of real observations and  numerical experimental designs, to inform the geometry of the numerical model. This information should be used to define model approximations that enable prior constraints to be defined. Inspired by the framework of this paper, such global approximations could be based on Gaussian mixtures of linear models estimated by regression (regression trees being a special form of such meta-models)\cite{Deleforge2015,Perthame2018,Li2018}, and benefit from, among other things, improvements in statistical estimation from shrinkage techniques \cite{Ledoit2004,Chen2010aa}. Indeed, these mixtures are becoming increasingly competitive computationally compared with kernel-based approaches, while avoiding their drawbacks (kernel ad-hoc choice, impossibility of inverting learned mappings \cite{Lawrence2005}). This aligns with the popularity of linear surrogates to approximate complex models or phenomena for the sake of examining their interpretability \cite{Hall2018}. In this context, the meaningfulness idea proposed in this article could potentially be adapted to constrain the model error to be informatively weaker than a stochastic input signal impacting the output signal. This should increase the relevance of the interpretability provided by these linear surrogates. 

In this sense, and beyond the benefits described in Proposition \ref{jeffreys},  we suggest that the idea proposed in this paper could inspire the improvement of Bayesian modeling choices in meta-modeling procedures, a common UQ framework. One problem that remains open is the selection of  covariance kernel structures when $g$ is replaced by a kriging Gaussian process $\hat{g}$, to which is added a second (zero-mean) Gaussian process as a prior on the so-called \emph{discrepancy} of $\hat{g}$ with respect to $g$ \cite{Kennedy2001,Xie2020}.  It is intuitive, for identifiability reasons, that  the additional discrepancy  should provide more contrast with respect to $\hat{g}$ than $\hat{g}$ with $g$ \cite{Brynjarsdottir2014}.  This can be interpreted as a reasonable prior constraint on the respective natures of the two covariance kernels. Nonetheless, to our knowledge no formal rule yet exists for this elicitation, but only rules of thumb based on the experience of researchers. As it proposes to define prior constraints from considerations on information degradation, the concept of  meaningfulness could be adapted to address this issue.

\section*{Acknowledgements}

The authors express their grateful thanks to Prs. Fabrice Gamboa and Jean-Michel Loubes (Institut de Math\'ematique de Toulouse)  for their  advices that helped to conduct this research. Drs Bertrand Iooss and Lorenzo Audibert (EDF Lab) are also thanked for fruitful discussions. Finally, the authors would like to thank the Associate Editors and two anonymous reviewers for their deep proofreading and precise comments that considerably helped to improve the tone and clarity of this article.

\appendix

\section{Grouped Sobol's indices}
\label{sobol}

\textcolor{black}{Closed Sobol's indices  \cite{sobol2001global}, based on Hoeffding's decomposition and the independence of all input variables, are one the most common indicators used in {global sensitivity analysis} of computer models. In the present situation we are interested by ranking the importance of two groups of input variables $X$ and $\varepsilon$ in the variance of output $Y^*=\hat{g}(Z)$ where $Z=(X,\varepsilon))\in \R^p\times \R^q$ and  $\hat{g}(Z)=g(X)+\varepsilon$. Each group of variables can contain dependent components but we assume that $X$ and $\varepsilon$ are independent. In this case, when $Y^*$ is univariate, so-called first-order \textit{grouped} Sobol's indices, defined by \cite{Jacques2006} and later studied by \cite{Broto2019}, generalize usual Sobol's indices: 
\begin{eqnarray}
S_X =  \frac{\rm{Var}[\E[Y^*|X]]}{\rm{Var}[Y^*]} & \text{and} & S_{\varepsilon} =  \frac{\rm{Var}[\E[Y^*|\epsilon]]}{\rm{Var}[Y^*]}. \label{grouped.sobol}
\end{eqnarray}
Other generalizations focus on the marginal output effects of each correlated input component in $Z$ (e.g., \cite{Mara2012}).  
The notion of {\it first order} is related to the fact that $S_X$ and  $S_{\varepsilon}$ do not measure the (second-order) impact of the interactions between $X$ and $\varepsilon$ to the output $Y^*$ \cite{Broto2019}. }

\textcolor{black}{More generally, when $Y^*$ is multidimensional, indices (\ref{grouped.sobol}) were extended by \cite{Gamboa2013} and \cite{GJKL14} under the strong hypothesis that the components of $Z$ are independent. }  
\textcolor{black}{Hoeffding's decomposition of $\tilde{g}$, assuming this independence between the  $Z_i$, leads to write an expansion of the following form \cite{sobol1993toto}
\begin{eqnarray*}
\tilde{g}(Z) & = & c + \tilde{g}_{\bf u}(Z_{\bf u}) + \tilde{g}_{/\bf u}(Z_{/\bf u}) +  \tilde{g}_{\bf u, /\bf u}(Z_{\bf u}, Z_{/\bf u})
\end{eqnarray*}
where ${\bf u}$ is any non-empty subset of ${\bf d}=\{1,\ldots,p+q\}$, $Z_{\bf u}  =  (Z_i, i \in {\bf u})$ and $Z_{/\bf u}  =  (Z_i, i \in {\bf d}/{\bf u})$, $c  =  \E[Y^*]$ and 
\begin{eqnarray*}
\tilde{g}_{\bf u} & = & \E[Y^*|Z_{\bf u}]-c, \\
\tilde{g}_{/\bf u} & = & \E[Y^*|Z_{/\bf u}]-c, \\
\tilde{g}_{\bf u, /\bf u} & = & Y^*-\tilde{g}_{\bf u}-\tilde{g}_{/\bf u}-c. 
\end{eqnarray*}
Thanks to $L^2$-orthogonality, 
\begin{eqnarray}
{\Cov[Y^*]} & = & {\Cov[\tilde{g}_{\bf u}(Z_{\bf u})]} + {\Cov[\tilde{g}_{/\bf u}(Z_{/\bf u})]} + {\Cov[\tilde{g}_{\bf u, /\bf u}(Z_{\bf u}, Z_{/\bf u})]}, \label{orthogo}
\end{eqnarray}
then the \emph{multivariate first-order grouped Sobol's indices} \cite{GJKL14} are defined by
\setcounter{equation}{2}
\begin{eqnarray}
S_{\bf u} & = & \frac{{\rm{Tr}}\left(I_q \Cov[\tilde{g}_{\bf u}(Z_{\bf u})]\right)}{{\rm{Tr}}\left(I_q \Cov[Y^*]\right)}.  \label{sobol.dim2}
\end{eqnarray}
In the specific case of model (\ref{inversion.problem}-\ref{eq:model_toy}), under the assumption of intern dependencies among the $X_i$ but since $g(X)$ is independent on $\varepsilon$, then
\begin{eqnarray}
{\Cov[Y^*]} & = & {\Cov[g(X)]} + {\Cov[\varepsilon]}. \label{a4}
\end{eqnarray}
Then Equation (\ref{orthogo}) still holds 
with $\tilde{g}_{\bf u}(Z_{\bf u})=g(X)$, $\tilde{g}_{/\bf u}(Z_{/\bf u})=\varepsilon$ and $\tilde{g}_{\bf u, /\bf u}(Z_{\bf u}, Z_{/\bf u})=0$, then 
the use of indices (\ref{sobol.dim2}) remains appropriate, their meaning being preserved. 
}

\section{The Fisher information}
\label{fisher}

The Fisher information $I_X(\theta)$ is a measure of the amount of information that an observable random variable $X$ carries about an unknown parameter $\theta \in \mathbb{R}$ upon which the probability of $X$ depends. This key statistical concept can help to quantify the uncertainty of a model or alternatively the amount of information carried by a model \cite{Ly2017}. \textcolor{black}{Especially, let us recall that $I_X(\theta)$ can be explained by a local differentiation of Shannon's entropy in the space of probability distributions, through the De Bruijn's identity \cite{Stam1959}.} If there exists no minimal set of necessary regularity conditions \emph{per se} for the existence of $I_X(\theta)$, most authors agree on the following sufficient conditions of existence, positivity and continuity in an subset of $\Theta$ (see for instance \cite{Hajek1972}, $\S$ 3.4): let $f(X; \theta)$ be the  probability density/mass of  $X$ conditional on the value of $\theta$. This function must be absolutely continuous in $\theta$ and the derivative ${\partial}f(x;\theta)/{\partial \theta} $ must exist for almost all $x$. \\ 

Then 
$$I_X(\theta)=\E\left[\left( \dfrac{\partial}{\partial \theta}\ln f(X, \theta)\right) ^2 \right] = \int  \left( \dfrac{\partial}{\partial \theta}\ln f(X, \theta)\right)^2 f(x,\theta)dx$$
or alternatively, if $\ln f(x, \theta)$ is twice differentiable with respect to $\theta$ and under some regularity conditions, 
$$I_X(\theta)=-\E\left[\dfrac{\partial^2 \ln f(X, \theta)}{\partial \theta^2} \right]. $$

Similarly, if $\theta=(\theta_1,...,\theta_N) \in \mathbb{R}^N$, the Fisher Information matrix $I_X(\theta)\in \mathbb{R}^{N\times N}$ is defined as 
$$\left( I_X(\theta)\right) _{ij}=-\E\left[\dfrac{\partial^2 \ln f(X, \theta)}{\partial \theta_i \theta_j}  \right]. $$

In the Gaussian case, where $X \sim \mathcal{N}_p\left( \mu(\theta), \Gamma(\theta)\right) $ and $\theta=(\theta_1,...,\theta_N)$, $\mu(\theta)= \left( \mu_1(\theta),...,\mu_p(\theta)\right) $, $\Sigma(\theta)= \left( \Gamma_{ij}(\theta)\right)_{1\leq i,j \leq p} $ , 
\begin{eqnarray}
\label{eq:FI_gaussian}
\left( I_X(\theta)\right)_{ij}= \dfrac{\partial \mu(\theta)^{T}}{\partial \theta_i}\Gamma(\theta)^{-1}\dfrac{\partial \mu(\theta)}{\partial \theta_j}+\dfrac{1}{2}\textrm{Tr}\left(\Gamma(\theta)^{-1}\dfrac{\partial \Gamma(\theta)}{\partial \theta_i}\Gamma(\theta)^{-1}\dfrac{\partial \Gamma(\theta)}{\partial \theta_j} \right) 
\end{eqnarray}
for all $1\leq i,j \leq N$.

If the mean and the covariance depend on two different parameters $\alpha \in \mathbb{R}^{l}$ and $\beta \in\mathbb{R}^{m}$, i.e $\theta:=(\alpha,\beta)$ and $X \sim \mathcal{N}_p(\mu(\alpha), \Gamma(\beta))$, then, using the Slepian-Bangs (SB) compact formula for the Fisher information matrix of a multivariate Gaussian distribution \cite{Kay1993} (Appendix 3C), 
\begin{eqnarray*}
I_X(\theta) & = & diag(I^{(1)}_{\alpha,\beta}, I^{(2)}_{\beta})
\end{eqnarray*}
 where
 \begin{eqnarray*}
\left( I^{(1)}_{\alpha,\beta}\right) _{1\leq i,j \leq l}=\dfrac{\partial \mu(\alpha)^{T}}{\partial \alpha_i}\Sigma(\beta)^{-1}\dfrac{\partial \mu(\alpha)}{\partial \alpha_j}
\end{eqnarray*}
and 
\setcounter{equation}{1}
\begin{eqnarray}
\left( I^{(2)}_{\beta} \right) _{1\leq i,j \leq m} & = & \dfrac{1}{2}\textrm{Tr}\left(\Sigma(\beta)^{-1}\dfrac{\partial \Sigma(\beta)}{\partial \beta_i}\Sigma(\beta)^{-1}\dfrac{\partial \Sigma(\beta)}{\partial \beta_j} \right). \label{fisher.beta}
\end{eqnarray}


\section{Proofs}


\begin{proof}[Proof of  Proposition \ref{prop1_bis}]
From Definition \ref{soboldef}, Equation (\ref{a4}) and for the model considered in (\ref{inversion.problem}-\ref{eq:model_toy}), 
\begin{eqnarray*}
S_X  =  \dfrac{
{\rm{Tr}}
\left(
\mbox{Cov}\left[Y\right]
\right)
}
{
{\rm{Tr}}
\left(
\mbox{Cov}\left[Y^*\right]
\right)
}, & & 
S_{\varepsilon}  =  \dfrac{{\rm{Tr}}\left(\mbox{Cov}\left[\varepsilon\right]\right)}{{\rm{Tr}}\left(\mbox{Cov}[Y^*]\right)}.
\end{eqnarray*}
Hence (\ref{sobol1}) is equivalent to ${\rm{Tr}}(a\Gamma a^T)> \Sigma$, Furthermore, the entropy of a multivariate normal distribution $X\sim\mathcal{N}(\mu, \Gamma)$ is 
\setcounter{equation}{0}
\begin{eqnarray}
\mathcal{E}(X) & = & \log ( \sqrt{( 2\pi e)\mid \Gamma \mid}). \label{eq:entropy_gaussian}
\end{eqnarray}
Then the proof is  straightforward. 
\end{proof}

\begin{proof}[Proof of  Proposition \ref{jeffreys}]
Jeffreys's prior for  the  ${\cal{N}}_p(\mu,\Gamma)$ distribution is, from \cite{yang1998},
\begin{eqnarray*}
\pi(\mu,\Gamma) & \propto & |\Gamma|^{-\frac{p+2}{2}}\1_{\{\Gamma\in S^{++}_p(\R)\}}\1_{\{\mu\in \R^p\}}
\end{eqnarray*}
and is such that $\int_{S^{++}_p} \pi(\Gamma) \ d\Gamma = \infty$. Obviously, restricting $\mu$ to $A$ implies $\pi(\Gamma)=\int_A \pi(\mu,\Gamma) \ d\mu < \infty$ for any finite $\Gamma$. Furthermore, rewrite $\Gamma=VDV^T$ where $D$ is diagonal and $V$ is the orthogonal matrix of eigenvectors. Then 
\begin{eqnarray*}
\int_{\Omega} \pi(\Gamma) \ d\Gamma & = & \left[ \int dV \left(\int_{\Omega_D} |D|^{-\frac{p+2}{2}} \ dD\right)\right]^{-1} 
\end{eqnarray*}
where $\Omega_D$ is a set of diagonal, positive matrices with upperly bounded values and such that $|D|>|\Sigma|/|a|^2$. Calculations conducted in Appendix B of \cite{Fu2015} show that $\int_{\Omega} \pi(\Gamma) \ d\Gamma<\infty$. Hence, the restriction of $\pi(\mu,\Gamma)$ over $A\times\Omega$ is proper. 
\end{proof}

\begin{proof}[Proof of  Proposition \ref{prop.gruff}]
Consider first the cross covariance
\begin{eqnarray}
\rm{Cov}\left(\tilde{g}_{x_0}(X),h_{x_0}(X) \right) & = & \rm{Cov}\left(\tilde{g}_{x_0}(X),g(X) - \tilde{g}_{x_0}(X)\right), \nonumber \\
& = & \E\left[\tilde{g}_{x_0}(X)g^T(X)\right] - \E\left[\tilde{g}_{x_0}(X)\tilde{g}^T_{x_0}(X)\right]- \E\left[\tilde{g}_{x_0}(X)\right]\E\left[{g}(X)\right]^T \nonumber  \\
& & \ \ \ - \E\left[\tilde{g}_{x_0}(X)\right]\E\left[\tilde{g}_{x_0}(X)\right]^T, \nonumber \\ \setcounter{equation}{2}
& = & \rm{Cov}\left(\tilde{g}_{x_0}(X),g(X) \right) - \rm{Cov}(\tilde{g}_{x_0}(X)). \label{crosscov}
\end{eqnarray}
Since, under H0, 
\begin{eqnarray*}
\rm{Cov}\left({g}(X) \right) & = & \rm{Cov}\left(\tilde{g}_{x_0}(X) \right) + \rm{Cov}\left(h_{x_0}(X) \right) + \rm{Cov}\left(\tilde{g}_{x_0}(X),h_{x_0}(X) \right) \\
& & \ \ \ + \rm{Cov}\left(\tilde{g}_{x_0}(X),h_{x_0}(X) \right)^T,
\end{eqnarray*}
then, from (\ref{crosscov}), 
\setcounter{equation}{4}
\begin{eqnarray}
{\rm{Tr}}\left[\rm{Cov}\left({g}(X) \right)\right] & = & {\rm{Tr}}\left[\rm{Cov}\left(\tilde{g}_{x_0}(X) \right)\right] + {\rm{Tr}}\left[\rm{Cov}\left(h_{x_0}(X) \right)\right] + 2{\rm{Tr}}\left[\rm{Cov}\left(\tilde{g}_{x_0}(X),h_{x_0}(X) \right)\right], \nonumber \\
& = &  {\rm{Tr}}\left[\rm{Cov}\left(h_{x_0}(X) \right)\right]  + 2{\rm{Tr}}\left[\rm{Cov}\left(\tilde{g}_{x_0}(X),g(X) \right)\right] \nonumber - {\rm{Tr}}\left[\rm{Cov}\left(\tilde{g}_{x_0}(X) \right)\right], \nonumber \\
& =&  {\rm{Tr}}\left[\rm{Cov}\left(h_{x_0}(X) \right)\right]  + 2{\rm{Tr}}\left[D_{g_{x_0}}\rm{Cov}\left(X,g(X) \right)\right] - {\rm{Tr}}\left[D_{g_{x_0}}\Gamma D^T_{g_{x_0}} \right] \nonumber \label{part1}
\end{eqnarray}
using cross covariance properties. When $q=1$, under H0 the multivariate version of Stein's identity (e.g., Lemma 1 in \cite{LIU1994247}) can be applied:
\begin{eqnarray}
\rm{Cov}\left(X,g(X) \right)  & = &  \Gamma \E\left[D_g(X)\right].  \label{stein.lemma}
\end{eqnarray}
Hence
\begin{eqnarray*}
{\rm{Tr}}\left[\rm{Cov}\left({g}(X) \right)\right] & = & {\rm{Tr}}\left[\rm{Cov}\left(h_{x_0}(X) \right)\right]  +{\rm{Tr}}\left[D_{g_{x_0}}\Gamma\left(2\E\left[D_g(X)\right] - D^T_{g_{x_0}}\right) \right].
\end{eqnarray*}
Since  ${\rm{Tr}}\left[\rm{Cov}\left(h_{x_0}(X) \right)\right]>0$, the proposed condition (\ref{res.propa})
implies 
\begin{eqnarray*}
{\rm{Tr}}\left[\mbox{Cov}(g(X))\right]& \geq & {\rm{Tr}}\left[\Sigma\right] \label{final.res1}
\end{eqnarray*}
which is similar to (\ref{sobol1}). 

\paragraph*{Remark.} Since, to our knowledge, there is today no multidimensional extension of Stein's lemma involving the Jacobian of $g$, providing a similar condition to (\ref{stein.lemma}) when $q>1$ seems an open problem. Alternative approaches to provide (and possibly refine) sufficient conditions could be based on trace inequalities \cite{Neumann1937,Liu2009}), possibly requiring that matrices involved in (\ref{res.propa}) are Hermitian (e.g., Ruhe's inequalities ; see \cite{Marshall2011}, p. 340-341). 
\end{proof}

\begin{proof}[Proof of  Theorem \ref{theo:fisher-info-gen}]
From Appendix \ref{fisher}, it can be seen that the Fisher matrix for (\ref{gen.model}) is block diagonal:
\begin{eqnarray*}
I_{{Y}_{\alpha}}(\theta) & = & \left(
\begin{array}{ll}
I_{{Y}_{\alpha}}(\mu) & {\bf 0} \\
{\bf 0} & I_{{Y}_{\alpha}}(\Gamma)
\end{array}\right)
\end{eqnarray*}
with 
\begin{eqnarray*}
\text{\bf (a)} \ \ \ I_{{Y}_{\alpha}}(\mu) & = & - \E_{{Y}_{\alpha}}\left[ \frac{\partial^2}{\partial \mu^2} \log f({Y}_{\alpha}|\theta)\right], \nonumber \\
                         & = & H^T \left(H\Gamma H^T + \alpha\Sigma\right)^{-1} H, \label{IYmu}\\
                         & =& H^T \Sigma^{-1/2}\left(\alpha I_q + \Psi_{\alpha}(\Gamma)\right)^{-1} \Sigma^{-1/2} H \ \ \text{if $\alpha=1$,} \\
                         &= & H^T \left(H\Gamma H^T\right)^{-1} H \hspace{2.9cm} \text{if $\alpha=0$.} \\
                         & & \hspace{-1cm}  \text{With  $\alpha\in\{0,1\}$ and   $ \Sigma^{0}=I_q$, a general expression is:}        \\   
                         & = & H^T\Sigma^{-\alpha/2} \left(\alpha I_q + \Psi_{\alpha}(\Gamma)\right)^{-1} \Sigma^{-\alpha/2}  H. \nonumber \\
\text{\bf (b)} \ \ \ I_{{Y}_{\alpha}}(\Gamma) & = & - \E_{{Y}_{\alpha}}\left[ \frac{\partial^2}{\partial \Gamma^2} \log f({Y}_{\alpha}|\theta)\right]. \nonumber
\end{eqnarray*}
Focusing on the second equation, notice that in the present case  
\begin{eqnarray*}
I_{{Y}_{\alpha}}(\Gamma) & = &   \dfrac{\partial^{2}}{\partial \Gamma^2}\mathbb{E}_{{Y}_{\alpha}}\left[- \log f({Y}_{\alpha}|\theta) \right] \ = \  \dfrac{\partial^{2}}{\partial \Gamma^2}  \mathcal{E}({Y}_{\alpha}) \\
& & \text{where  $\mathcal{E}({Y}_{\alpha})$ is the entropy (see Proposition \ref{prop1_bis})}, \\
& = &   -\dfrac{1}{2}\dfrac{\partial^{2}}{\partial \Gamma^2}\log\left| H\Gamma H^T + \alpha\Sigma\right|, \\
& = &  -\dfrac{1}{2}\dfrac{\partial^{2}}{\partial \Gamma^2} \log \left|  \Psi_{\alpha}(\Gamma) + \alpha I_q\right|.
\end{eqnarray*}
Since $\Psi_{\alpha}(\Gamma)$ is a real symmetrical matrix in $\R^{q,q}$, it is diagonalisable. Then real eigenvalues $(\lambda^{\Psi_{\alpha}(\Gamma)}_i)_{1\leq i \leq q}$ exist, which are the solutions of the characteristic polynomial of  $-\Psi_{\alpha}(\Gamma)$:
\begin{eqnarray*}
p_{\Gamma}(x) & = & \left| -\Psi_{\alpha}(\Gamma) + \alpha xI_q\right|, \\
& = &  \prod\limits_{i=1}^q \left(\lambda^{\Psi_{\alpha}(\Gamma)}_i + \alpha x\right).
\end{eqnarray*}
Hence 
\begin{eqnarray*}
I_{{Y}_{\alpha}}(\Gamma) & = & -\dfrac{1}{2}\dfrac{\partial^{2}}{\partial \Gamma^2}  p_{\Gamma}(1), \\
& = & -\dfrac{1}{2}\sum\limits_{i=1}^q \dfrac{\partial^{2}}{\partial \Gamma^2} \log\left(\alpha + \lambda^{\Psi_{\alpha}(\Gamma)}_i\right).
\end{eqnarray*}
Note that since $\Psi_{\alpha}(\Gamma)$ is symmetric, Theorem 1 in \cite{Magnus1985} ensures the existence of the second order derivative of the $\lambda^{\Psi_{\alpha}(\Gamma)}_i$ (or some differentiable function of it) with respect to $\Gamma$, in a neighborhood of $\Gamma$. Hence the last equation makes sense.
\end{proof}
\begin{proof}[Proof of Proposition \ref{going.further}]
Considering the situation $\Gamma=diag(\tau^2_1,\ldots,\tau^2_p)$ of dimension $p$, first notice the following result.
With $H=[h_{ij}]_{i,j}$ for $1\leq i\leq q$ and $1\leq j\leq p$, the matrix $\Psi_{\alpha}(\Gamma)$ becomes
\begin{eqnarray}
\Psi_{\alpha}(\tau^2_1,\ldots,\tau^2_p)  =   \Sigma^{-\alpha/2} \left[\sum\limits_{k=1}^p \tau^2_k h_{ik}h_{jk}\right]_{i,j} \Sigma^{-\alpha/2} 
 =   \sum\limits_{k=1}^{p} \tau^2_k \tilde{A}_{k,\alpha}. \label{general.tauk.formula}
\end{eqnarray}
This provides Expression (\ref{fisher.term.simp.1}) using Theorem \ref{theo:fisher-info-gen}.
Besides, the former (crude) notation $\partial / \partial \Gamma$ in Theorem \ref{theo:fisher-info-gen} becomes $(\partial / \partial \tau^2_i)_{1\leq i \leq p}$. The $(k,l)-$component of the Fisher matrix $I_{{Y}_{\alpha}}(\Gamma)$  thus becomes
\begin{eqnarray*}
I_{{Y}_{\alpha}}(\Gamma)_{k,l} & = & -\dfrac{1}{2}\sum\limits_{i=1}^q \dfrac{\partial}{\partial \tau^2_k}\dfrac{\partial}{\partial \tau^2_l} \log\left(\alpha + \lambda^{\Psi_{\alpha}(\Gamma)}_i\right), \\
& = & -\dfrac{1}{2}\sum\limits_{i=1}^q \dfrac{\partial}{\partial \tau^2_k} 
      \left\{
      \left(\dfrac{\partial}{\partial \tau^2_l}\lambda_i^{\Psi_{\alpha}(\Gamma)}\right) \frac{1}{\left(\alpha + \lambda_i^{\Psi_{\alpha}(\Gamma)}\right)}\right\}, \\
& = & -\dfrac{1}{2}\sum\limits_{i=1}^q 
\left\{ 
 \left(\dfrac{\partial^2}{\partial \tau^2_l \partial \tau^2_k}\lambda_i^{\Psi_{\alpha}(\Gamma)}\right) \frac{1}{\left(\alpha + \lambda_i^{\Psi_{\alpha}(\Gamma)}\right)} \right.\\
 & & \ \ \ - \ 
 \left. 
 \left(\dfrac{\partial}{\partial \tau^2_l}\lambda_i^{\Psi_{\alpha}(\Gamma)}\right) \left(\dfrac{\partial}{\partial \tau^2_k}\lambda_i^{\Psi_{\alpha}(\Gamma)}\right) \frac{1}{\left(\alpha + \lambda_i^{\Psi_{\alpha}(\Gamma)}\right)^2} 
\right\}.
\end{eqnarray*}
The statement (\ref{fisher.term.simp.2}) follows straightforwardly from Lemma \ref{derivative.lambda.2}.
\end{proof}

\begin{lemma}\label{derivative.lambda.2}
Using the notations defined in the proof of Proposition \ref{going.further}, one has
\begin{eqnarray}
\frac{\partial}{\partial \tau^2_k} \lambda^{\Psi_{\alpha}(\Gamma)}_i  & =  & v^T_{\alpha,i} \tilde{A}_{k,\alpha} v_{\alpha,i}, \label{control.derivative.1} \\
\frac{\partial^2}{\partial \tau^2_k \partial \tau^2_l} \lambda^{\Psi_{\alpha}(\Gamma)}_i &= & v^T_{\alpha,i} \tilde{A}_{l,\alpha} \Psi^+_i(\Gamma)  \tilde{A}_{k,\alpha} v_{\alpha,i}  +  v^T_{\alpha,i}  \tilde{A}_{k,\alpha}  \Psi^+_i(\Gamma)  \tilde{A}_{l,\alpha} v_{\alpha,i}
\end{eqnarray}
where $A^+$ is the generalized Moore-Penrose pseudo-inverse of $A$ and 
\begin{eqnarray*}
\Psi_i(\Gamma) & = &  \lambda^{\Psi_{\alpha}(\Gamma)}_i I_q- \Psi_{\alpha}(\Gamma).
\end{eqnarray*}
\end{lemma}

\begin{proof}[Proof of Lemma \ref{derivative.lambda.2}]
By definition, for $i\in\{1,\ldots,q\}$,
\begin{eqnarray}
\Psi_{\alpha}(\Gamma) v_{\alpha,i} & = & \lambda^{\Psi_{\alpha}(\Gamma)}_i v_{\alpha,i}.\label{eigenvalues.def}
\end{eqnarray}
Hence, for $k\in\{1,\ldots,p\}$, 
\begin{eqnarray*}
\frac{\partial}{\partial \tau^2_k} \left(\Psi_{\alpha}(\Gamma) v_{\alpha,i} \right) & = & \frac{\partial}{\partial \tau^2_k} \left(\Psi_{\alpha}(\Gamma) \right) v_{\alpha,i}  +  \Psi_{\alpha}(\Gamma)  \frac{\partial}{\partial \tau^2_k}  v_{\alpha,i}, \\
& = & \tilde{A}_{k,\alpha} v_{\alpha,i} + \Psi_{\alpha}(\Gamma)  \frac{\partial}{\partial \tau^2_k}  v_{\alpha,i} \ \ \text{from (\ref{general.tauk.formula}),} \\
& = & \frac{\partial}{\partial \tau^2_k} \left(\lambda^{\Psi_{\alpha}(\Gamma)}_i v_{\alpha,i}\right) \ \ \text{from (\ref{eigenvalues.def}),} \\
& = & \left(\frac{\partial}{\partial \tau^2_k} \lambda^{\Psi_{\alpha}(\Gamma)}_i\right) v_{\alpha,i} + \lambda^{\Psi_{\alpha}(\Gamma)}_i \frac{\partial}{\partial \tau^2_k}  v_{\alpha,i}.
\end{eqnarray*}
Then, left-multiplying by $v^T_{\alpha,i}$ the above expressions and noticing that $v^T_{\alpha,i}v_{\alpha,i}=1$, 
\begin{eqnarray*}
\frac{\partial}{\partial \tau^2_k} \lambda^{\Psi_{\alpha}(\Gamma)}_i &= & v^T_{\alpha,i} \tilde{A}_{k,\alpha} v_{\alpha,i} + v^T_{\alpha,i}\Psi_{\alpha}(\Gamma)  \frac{\partial}{\partial \tau^2_k}  v_{\alpha,i}  - \lambda^{\Psi_{\alpha}(\Gamma)}_i  v^T_{\alpha,i} \frac{\partial}{\partial \tau^2_k}  v_{\alpha,i}.
\end{eqnarray*}
Furthermore, from (\ref{eigenvalues.def}), $v^T_{\alpha,i}\Psi_{\alpha}(\Gamma)=\lambda^{\Psi_{\alpha}(\Gamma)}_i v^T_{\alpha,i}$. Hence
\begin{eqnarray*}
\frac{\partial}{\partial \tau^2_k} \lambda^{\Psi_{\alpha}(\Gamma)}_i &= & v^T_{\alpha,i} \tilde{A}_{k,\alpha} v_{\alpha,i}.
\end{eqnarray*}
The result (\ref{control.derivative.1}) is then straightforward. 
Besides, from \cite{Magnus1985}, 
\begin{eqnarray}
\frac{\partial}{\partial \tau^2_k} v_{\alpha,i} & = & \left\{ \lambda^{\Psi_{\alpha}(\Gamma)}_i I_q- \Psi_{\alpha}(\Gamma)\right\}^+ \frac{\partial}{\partial \tau^2_k} \left(\Psi_{\alpha}(\Gamma) \right) v_{\alpha,i}, \nonumber \\
& = & \Psi^+_i(\Gamma) \tilde{A}_{k,\alpha} v_{\alpha,i}.\label{deriv.eigenvec}
\end{eqnarray}
Then it comes
\begin{eqnarray*}
\frac{\partial^2}{\partial \tau^2_k \partial \tau^2_l} \lambda^{\Psi_{\alpha}(\Gamma)}_i &= & \frac{\partial}{\partial \tau^2_l} \left(v^T_{\alpha,i} \tilde{A}_{k,\alpha} v_{\alpha,i}\right), \\
& = & \left( \frac{\partial}{\partial \tau^2_l} v^T_{\alpha,i}\right)  \tilde{A}_{k,\alpha} v_{\alpha,i}  + v^T_{\alpha,i}  \tilde{A}_{k,\alpha} \left( \frac{\partial}{\partial \tau^2_l} v_{\alpha,i}\right)
\end{eqnarray*}
since $\frac{\partial}{\partial \tau^2_l}  \tilde{A}_{k,\alpha} = 0$. Hence, from (\ref{deriv.eigenvec}), since $\Psi^+_i(\Gamma)=\left\{ \lambda^{\Psi_{\alpha}(\Gamma)}_i I_q- \Psi_{\alpha}(\Gamma)\right\}^+$ and $\tilde{A}_{l,\alpha}$ are symmetric, 
\begin{eqnarray*}
\frac{\partial^2}{\partial \tau^2_k \partial \tau^2_l} \lambda^{\Psi_{\alpha}(\Gamma)}_i &= & v^T_{\alpha,i} \tilde{A}_{l,\alpha} \Psi^+_i(\Gamma)  \tilde{A}_{k,\alpha} v_{\alpha,i}  +  v^T_{\alpha,i}  \tilde{A}_{k,\alpha}  \Psi^+_i(\Gamma) \tilde{A}_{l,\alpha} v_{\alpha,i}.
\end{eqnarray*} 
Furthermore,  $(\lambda^{\Psi_{\alpha}(\Gamma)}_i -\lambda^{\Psi_{\alpha}(\Gamma)}_j)_j$ is the eigenvalue of  $\Psi_i(\Gamma)$ related to $v_{\alpha,j}$. Since $\Psi^+_i(\Gamma)$ is the matrix for which $\Psi^+_i(\Gamma)\Psi_i(\Gamma) x = x$ for all $x$ in the row space of $\Psi_i(\Gamma)$, then for all $j$ such that  $\lambda^{\Psi_{\alpha}(\Gamma)}_j \neq \lambda^{\Psi_{\alpha}(\Gamma)}_i$ it comes
\begin{eqnarray*}
\Psi^+_i(\Gamma) v_{\alpha,j} & = & \frac{1}{\lambda^{\Psi_{\alpha}(\Gamma)}_i -\lambda^{\Psi_{\alpha}(\Gamma)}_j} \Psi^+_i(\Gamma)\Psi_i(\Gamma) v_{\alpha,j} \ = \  \frac{1}{\lambda^{\Psi_{\alpha}(\Gamma)}_i -\lambda^{\Psi_{\alpha}(\Gamma)}_j} v_{\alpha,j}
\end{eqnarray*}
hence the nonzero eigenvalues of $\Psi^+_i(\Gamma)$ are the $(\lambda^{\Psi_{\alpha}(\Gamma)}_i -\lambda^{\Psi_{\alpha}(\Gamma)}_j)^{-1}$ and $\Psi^+_i(\Gamma)$ has $n_i=\mbox{Card}(1\leq j \leq p, \ \lambda^{\Psi_{\alpha}(\Gamma)}_j=\lambda^{\Psi_{\alpha}(\Gamma)}_i)$ zero eigenvalues.
\end{proof}

\begin{proof}[Proof of  Corollary \ref{diagonalite.aa}]
Assume that for any $(k,l)\in\{1,\ldots,p\}^2$, $H_k H^T_l$ is symmetric. From Lemma \ref{symmetry.h}, then $H_k H^T_k$ and $H_l H^T_l$ commute. Assume furthermore that $\Sigma=diag(\sigma^2_1,\ldots,\sigma^2_q)$. Then $\Sigma^{-\alpha/2}=(\sigma^{-\alpha}_{1},\ldots,\sigma^{-\alpha}_{q})$ and for any $(k,l)\in\{1,\ldots,p\}^2$,
\begin{eqnarray*}
& \tilde{A}_{k,\alpha}\tilde{A}_{l,\alpha} - \tilde{A}_{l,\alpha}\tilde{A}_{k,\alpha}   =   \left(\sum\limits_{u=1}^q h_{u,k} h_{u,l} \sigma^{-2\alpha}_u\right) \left[\sigma^{-\alpha}_i\sigma^{-\alpha}_j \left(h_{i,k}h_{j,l}-h_{i,l}h_{j,k}\right)\right]_{1\leq i,j\leq q}, \\
&  \hspace{-5.3cm} =   0,
\end{eqnarray*}
hence $\tilde{A}_{k,\alpha}$ and $\tilde{A}_{l,\alpha}$ commute. In such a case, from \cite{Knutson2001} and (\ref{general.tauk.formula}), for $i\in\{1,\ldots,q\}$
\begin{eqnarray*}
\lambda^{\Psi_{\alpha}(\Gamma)}_i & = & \sum\limits_{j=1}^{p} \lambda^{\tilde{A}_{j,\alpha}}_i \tau^2_j. 
\end{eqnarray*}
(Note that a  general result on the relation between the eigenvalues of a sum with the sum of eigenvalues is provided by the Knutson-Tao theorem in \cite{Knutson1999}). Hence, from Theorem \ref{theo:fisher-info-gen}, $I_{{Y}_{\alpha}}(\Gamma)$ is the following $p\times p$ matrix
\begin{eqnarray*}
I_{{Y}_{\alpha}}(\Gamma) & = & -\dfrac{1}{2}\sum\limits_{i=1}^q\left[ \dfrac{\partial^{2}}{\partial \tau^2_k \partial \tau^2_l}   \log\left(\alpha + \sum\limits_{j=1}^{p} \lambda^{\tilde{A}_{j,\alpha}}_i \tau^2_j\right) \right]_{1 \leq k,l \leq p}.
\end{eqnarray*}
For any $j\in\{1,\ldots,p\}$, the $q$ eigenvalues $\lambda^{\tilde{A}_{j,\alpha}}_i$ of the symmetric matrix $\tilde{A}_{j,\alpha}$ are all 0 except 
$\lambda^{\tilde{A}_{j,\alpha}}_{{i^*_j}}  =  \lambda^{\tilde{A}_{j,\alpha}}_{\sup} = \| \tilde{A}_{j,\alpha} \|_2$ 
since $H_j H^T_j$ is of rank 1. Furthermore, since $\tilde{A}_{j,\alpha}=(\Sigma^{-\alpha/2} H_j)(\Sigma^{-\alpha/2} H_j)^T$, then
$ \| \tilde{A}_{j,\alpha} \|_2  = \| \Sigma^{-\alpha/2} H_j\|^2_2$. Hence
\begin{eqnarray*}
I_{{Y}_{\alpha}}(\Gamma) & = & -\dfrac{1}{2}\left[ \dfrac{\partial^{2}}{\partial \tau^2_k \partial \tau^2_l} \log\left(\alpha + \sum\limits_{j=1}^{p} \| \Sigma^{-\alpha/2} H_j\|^2_2 \tau^2_j\right) \right]_{1 \leq k,l \leq p}, \\
& = & \left[ \dfrac{1}{2} \frac{\| \Sigma^{-\alpha/2} H_k\|^2_2 \| \Sigma^{-\alpha/2} H_l\|^2_2 }{ \left(\alpha + \sum\limits_{j=1}^{p}\| \Sigma^{-\alpha/2} H_j\|^2_2 \tau^2_j \right)^2}\right]_{1 \leq k,l \leq p}.
\end{eqnarray*}
It is equivalent to write
\begin{eqnarray*}
I_{{Y}_{\alpha}}(\Gamma) & = & \frac{\Lambda_{{\alpha}} \Lambda^T_{{\alpha}} }{2\left(\alpha + {\rm{Tr}}\left(\tilde{\Lambda}_{\alpha} \Gamma\right) \right)^2} 
\end{eqnarray*}
where $\Lambda_{{\alpha}}=(\| \Sigma^{-\alpha/2} H_1\|^2_2,\ldots, \| \Sigma^{-\alpha/2} H_p\|^2_2)^T$ and  $\tilde{\Lambda}_{{\alpha}}=diag \Lambda_{{\alpha}}$. 
\end{proof}

 \begin{lemma}\label{symmetry.h}
Denote $H_k=(h_{i,k})_{1\leq i \leq q}$. Then the two following assumptions are equivalent:
\begin{description}
\item[(i)] For any $(k,l)\in\{1,\ldots,p\}^2$, $H_k H^T_l$ is symmetric;
\item[(ii)] For any $(k,l)\in\{1,\ldots,p\}^2$, $H_k H^T_k$ and $H_l H^T_l$ commute.
\end{description}
\end{lemma}

\begin{proof}[Proof of Lemma \ref{symmetry.h}]
For any $(k,l)\in\{1,\ldots,p\}^2$, notice that
\begin{eqnarray*}
\left(H_k H^T_k\right)\left(H_l H^T_l\right) & = & \left[\sum\limits_{u=1}^q h_{i,k} h_{u,k} h_{u,l} h_{j,l}\right]_{1\leq i,j \leq q}, \\
& = & \left(\sum\limits_{u=1}^q h_{u,k}  h_{u,l}\right)  \left[h_{i,k}h_{j,l}\right]_{1\leq i,j \leq q}
\end{eqnarray*}
and similarly
\begin{eqnarray*}
\left(H_l H^T_l\right)\left(H_k H^T_k\right)
& = & \left(\sum\limits_{u=1}^q h_{u,k}  h_{u,l}\right)  \left[h_{i,l}h_{j,k}\right]_{1\leq i,j \leq q}.
\end{eqnarray*}
Then the last two equations are equal if and only if 
\begin{eqnarray*}
H_k H^T_l & = & \left[h_{i,k}h_{j,l}\right]_{i,j} \ = \  \left[h_{i,l}h_{j,k}\right]_{i,j} \ = \ \left(H_k H^T_l\right)^T \ = \ H_l H^T_k.
\end{eqnarray*}
\end{proof}
\begin{proof}[Proof of  Proposition \ref{NC.A}]
Denote $A=\Psi_0(\Gamma)=H\Gamma H^T$ and $C=H\Gamma H^T + \Sigma$. With $(A,C,\Sigma)\in S_q^{++}(\R)$, then $A\prec C$.
For all $i\in\{1,\ldots,q\}$, one has $\lambda^{A^{-1}C}_i  =  \lambda^{I_q + A^{-1}\Sigma}_i \ > \ 0$. Hence, from Theorem 1 in \cite{Liski1996},  
$V  =  (A^{-1}-C^{-1})^{-1}$ is invertible and symmetric. Hence there exists $N\succeq 0$ in $\R^{q,q}$ such that $V^{-1}=N^TN$. Consequently, from Theorem \ref{theo:fisher-info-gen}, 
\begin{eqnarray*}
I_{Y_0}(\mu) - I_{Y_1}(\mu) & = &  H^T \left(A^{-1}-C^{-1}\right) H, \\
  & = & H^T V^{-1} H, \\
   & = & (NH)^T (NH) \ \succeq \ 0 \ \ \ \text{which proves (\ref{CN.1}).}
\end{eqnarray*}
\end{proof}

\begin{proof}[Proof of  Proposition \ref{NC.C}]
Denote $\Gamma=(\tau^2_1,\ldots,\tau^2_p)$. 
Since $\Psi_0(\Gamma)$ and $\Sigma^{1/2}$ are co-diagonalizable, and since they belong to $S^{++}_q(\R)$, there exists an invertible matrix $P$ and two invertible diagonal matrices $(\Lambda_{0},\Lambda_1)$ of rank $q$, containing the eigenvalues of $H\Gamma H^T$ and $\Sigma^{1/2}$, such that  $\Psi_0(\Gamma)=P\Lambda_0P ^{-1}$ and  $\Sigma^{1/2}=P\Lambda_1P ^{-1}$. Then 
\begin{eqnarray*}
\Psi_1(\Gamma) & = & \Sigma^{-1/2} \Psi_0(\Gamma)\Sigma^{-1/2} \ = \ P \Lambda^{-1}_1\Lambda_0 \Lambda^{-1}_1P^{-1}
\end{eqnarray*}
and for $i=1,\ldots,q$
\begin{eqnarray*}
\lambda^{\Psi_1(\Gamma)}_i & = & \lambda^{H\Gamma H^T}_i \left(\lambda^{\Sigma^{1/2}}\right)^{-2} \ = \ \lambda^{\Psi_0(\Gamma)}_i /  \lambda^{\Sigma}_i.
\end{eqnarray*}
In this situation, from Theorem \ref{theo:fisher-info-gen}, for $(k,l)\in\{1,\ldots,p\}$,
\begin{eqnarray*}
I_{{Y}_{\alpha}}(\Gamma)_{k,l} & = & \frac{1}{2} \sum\limits_{i=1}^q \omega^{(\alpha)}_i(\Gamma) \frac{\partial \lambda^{\Psi_0(\Gamma)}_i}{\partial \tau^2_k} \frac{\partial \lambda^{\Psi_0(\Gamma)}_i}{\partial \tau^2_l}, 
\end{eqnarray*} 
with $ \omega^{(\alpha)}_i(\Gamma) =   \left(\alpha\lambda^{\Sigma}_i + \lambda^{\Psi_0(\Gamma)}_i\right)^{-2}$. Given now $z=(z_1,\ldots,z_p)\in\R^p$, $z\neq {\bf 0}$, then
\begin{eqnarray*}
z^T \left(I_{{Y}_{0}}(\Gamma)-I_{{Y}_{1}}(\Gamma)\right) z & = & \frac{1}{2}\sum\limits_{i=1}^q \left\{\omega^{(0)}_i(\Gamma)-\omega^{(1)}_i(\Gamma)\right\} z^T M_i M^T_i z
\end{eqnarray*}
with $M_i=(\partial \lambda^{\Psi_0(\Gamma)}_i /{\partial \tau^2_1},\ldots,\partial \lambda^{\Psi_0(\Gamma)}_i/{\partial \tau^2_p})$. With $z^T M_i M^T_i z\geq 0$ and  $\omega^{(0)}_i(\Gamma) > \omega^{(1)}_i(\Gamma)$ $\forall i \in \{1,\ldots,q\}$, then
$z^T \left(I_{{Y}_{0}}(\Gamma)-I_{{Y}_{1}}(\Gamma)\right) z  \geq  0$ which proves $I_{{Y}_{0}}(\Gamma) \succeq I_{{Y}_{1}}(\Gamma)$.
\end{proof}
\begin{proof}[Proof of Proposition \ref{NC.B}]
Reusing the notations of the proof of Proposition \ref{NC.A}, denote $W_c  =  C^{-1}-\frac{1}{c} A^{-1}$. With $W^{-1}_c$ symmetric and real, there exists an orthogonal $q\times q-$matrix $Q$ and a diagonal matrix $\Lambda$ such that $W^{-1}_c=Q \Lambda Q^T$. Hence
\begin{eqnarray*}
I_{Y_1}(\mu) - \frac{1}{c}I_{Y_0}(\mu) & = & H^T W_c H, \\
& = & (H^TQ) \Lambda (H^TQ)^T \ \succeq \ 0
\end{eqnarray*}
if and only if each diagonal element of $\Lambda$ is positive or null. Hence $W_c\succeq 0$ is a NSC for  (\ref{CN.2}), or similarly
\begin{eqnarray}
C^{-1} & \succeq & \frac{1}{c} A^{-1}. \label{ineq.toto} 
\end{eqnarray}
Since, $\forall i \in\{1,\ldots,q\}$,
$\lambda_i^{cA C^{-1}}  =  \lambda_i^{c\left(I_q + \Sigma (H\Gamma H^T)^{-1}\right)}  \geq  0$,
then from Theorem 1 in \cite{Liski1996}, (\ref{ineq.toto}) is equivalent to 
\begin{eqnarray*}
cA & \succeq & C=A+\Sigma \ \ \Leftrightarrow \ \ H\Gamma H^T \ \succeq \ \frac{1}{c-1} \Sigma.
\end{eqnarray*}
Furthermore, since $\lambda^{\Gamma}_{\sup}  I_p\succeq \Gamma \succeq \lambda^{\Gamma}_{\inf} I_p$ \cite{Zhan2002},  $\lambda^{\Gamma}_{\sup} =\|\Gamma\|_2$, then
$D_1  =  \Gamma - \lambda^{\Gamma}_{\inf}  I_p \ \succeq \ 0$ 
and there exists a  matrix $\tilde{Q}\in\R^{p,p}$ such that $D_1=\tilde{Q}\tilde{Q}^T$. Hence $A=H\Gamma H^T=(H\tilde{Q})(H\tilde{Q})^T \succeq 0$. From the above expression, it follows that   $H\Gamma H^T \succeq \lambda^{\Gamma}_{\inf} HH^T$. Replacing $D_1$ by $D_2=\|\Gamma\|_2  I_p -\Gamma$ and following the same rationale, it comes $H\Gamma H^T \preceq \|\Gamma\|_2 HH^T$. Then a necessary (\textit{resp.} sufficient) condition for  (\ref{fisher.0.grosgris}) is
\begin{eqnarray*}
\|\Gamma\|_2 HH^T   \ \ \  \left(resp. \  \ \lambda^{\Gamma}_{\inf} HH^T\right) & \succeq & \frac{1}{c-1} \Sigma
\end{eqnarray*}
or equivalently, from Lemma 1 in \cite{Hauke1994}: 
\begin{eqnarray*}
\|\Gamma\|_2  \ \ \  \left(resp. \ \   \lambda^{\Gamma}_{\inf} \right)  & \geq & \frac{1}{c-1} \lambda_{\sup}^{(HH^T)^{-1}\Sigma}, \\
& & \ \ = \ \  \frac{1}{c-1}\| \left(HH^T\right)^{-1} \Sigma\|_2
\end{eqnarray*}
(note besides that $\lambda^{\Gamma}_{\inf} =\|\Gamma^{-1}\|^{-1}_2$ if $\Gamma\succ 0$).
\end{proof}
\begin{proof}[Proof of Proposition \ref{NC.D}]
From (\ref{I.mu.cov}), (\ref{CN.4}) is true if and only if
\begin{eqnarray*}
 \frac{\Lambda_{1} \Lambda^T_{1}}{\left(\alpha + {\rm{Tr}}\left(\tilde{\Lambda}_{1} \Gamma\right) \right)^2} & \succeq & \frac{1}{c} \frac{\Lambda_{0} \Lambda^T_{0}}{\left(\alpha + {\rm{Tr}}\left(\tilde{\Lambda}_{0} \Gamma\right) \right)^2}
\end{eqnarray*}
which, from Lemma 1 in \cite{Hauke1994} and since $\Lambda_{0} \Lambda^T_{0}\succeq 0$, is equivalent to 
\begin{eqnarray}
\sqrt{c} \frac{{\rm{Tr}}\left(\tilde{\Lambda}_0 \Gamma\right)}{1 + {\rm{Tr}}\left(\tilde{\Lambda}_1 \Gamma\right)}  \geq  \rho
&
\text{or similarly} &  
{\rm{Tr}}(U\Gamma)  \geq  \rho \label{cond.cur.rho}
\end{eqnarray}
where $\rho  =  \sqrt{\left\| \left(\Lambda_{1} \Lambda^T_{1}\right)^+  \left(\Lambda_{0} \Lambda^T_{0}\right) \right\|_2}$ and $ U   =   \sqrt{c}\tilde{\Delta}_0 - \rho \tilde{\Delta}_1$. Noticing that $U=diag(u_1,\ldots,u_p)$ where
\begin{eqnarray*}
u_i & = &  \sqrt{c}\left\| H_i \right\|^2_2 - \rho \left\| \Sigma^{-1/2} H_i \right\|^2_2 
\end{eqnarray*}
and $\Sigma^{-1/2} H_i \preceq \|\Sigma^{-1/2}\|_2 H_i$ since $\Sigma^{-1/2}$ is diagonal, then $u_i \geq \| H_i \|^2_2(\sqrt{c} - \rho \|\Sigma^{-1/2}\|^2_2)$ and consequently
\begin{eqnarray}
U & \succeq & \tilde{\Lambda}_0 \left(\sqrt{c} - \rho \|\Sigma^{-1/2}\|^2_2\right). \label{lower.bound.U}
\end{eqnarray}
Since $\Lambda_{1} \Lambda^T_{1}$ is diagonalizable and of rank 1,  there exists an orthogonal matrix $P$ such that $\Lambda_{1} \Lambda^T_{1}=P W P^{-1}$ where  $W=diag(0,\ldots,0,\|\Lambda_1\|^2_2, 0,\ldots,0)$ with $,\|\Lambda_1\|_2>0$. Then the Moore-Penrose pseudo-inverse $(\Lambda_{1} \Lambda^T_{1})^+= P\tilde{W} P^{-1}$ where $\tilde{W}=diag(0,\ldots,0,\|\Lambda_1\|^{-2}_2, 0,\ldots,0)$. Hence $\|(\Lambda_{1} \Lambda^T_{1})^+\|_2 = \|\Lambda_1\|^{-2}_2$. Consequently, 
\begin{eqnarray*}
\rho & \leq & \sqrt{\left\| \left(\Lambda_{1} \Lambda^T_{1}\right)^+\right\|_2 \left\| \left(\Lambda_{0} \Lambda^T_{0}\right)\right\|_2 } \ = \  \|\Lambda_1\|^{-1}_2 \|\Lambda_0\|_2. 
\end{eqnarray*}
Then, from (\ref{lower.bound.U}),
\begin{eqnarray*}
U & \succeq & \tilde{\Lambda}_0 \left(\sqrt{c} - \gamma \right)
\end{eqnarray*}
with
\begin{eqnarray*}
\gamma & = &  \|\Sigma^{-1/2}\|^2_2\|\Lambda_1\|^{-1}_2 \|\Lambda_0\|_2, \\
              & = & \left(\frac{\sum\limits_{j=1}^p  \|\Sigma^{-1/2}\|^4_2 \left\| H_i \right\|^4_2}{\sum\limits_{j=1}^p \left\|\Sigma^{-1/2} H_i \right\|^4_2} \right)^{1/2}.
\end{eqnarray*}
Then, if $c>\gamma^2$, then $U\succ 0$. Then, with $(U,\Gamma)$ diagonal, from the Cauchy-Schwarz inequality, ${\rm{Tr}}(U\Gamma)\leq \|\Gamma\|_2 {\rm{Tr}}(U) \leq  {\rm{Tr}}(U){\rm{Tr}}(\Gamma)$. Consequently, a 
 NC for (\ref{cond.cur.rho}) is
\begin{eqnarray*}
\|\Gamma\|_2  & \geq & \frac{\rho}{ {\rm{Tr}}\left(U\right)} \ = \ \frac{\rho}{\sqrt{c}\| \Lambda_0 \|^2_2 - \rho \| \Lambda_1 \|^2_2}.
\end{eqnarray*}
Furthermore, since $U\Gamma   \succeq   \tau^2_{\inf} U$, the corresponding SC for (\ref{cond.cur.rho}) is
\begin{eqnarray*}
 \tau^2_{\inf} & \geq & \frac{\rho}{ {\rm{Tr}}\left(U\right)}.
\end{eqnarray*}
\end{proof}
\begin{proof}[Proof of Proposition \ref{all.res}]
Theorem \ref{theo:fisher-info-gen} and Corollary \ref{diagonalite.1} show that a NC for \eqref{eq:fisher-condition2}, assuming $\Gamma=\tau^2 I_p$, 
is  
\begin{eqnarray}
\dfrac{q}{2}\left( \dfrac{1}{\tau^2}\right) ^2 >\dfrac{1}{2}\left( \dfrac{1}{\tau^2}\right) ^2\sum_{i=1}^{q}\left(\dfrac{\lambda_i^{\Psi_1}}{\lambda_i^{\Psi_1}+1/\tau^2} \right)^2 > \dfrac{1}{c}\dfrac{q}{2}\left( \dfrac{1}{\tau^2}\right) ^2. \label{NC.fishfish}
\end{eqnarray}
The left-side inequality of \eqref{NC.fishfish} is always satisfied since $\tau^2>0$ and $$\dfrac{\lambda_i^{\Psi_1}}{\lambda_i^{\Psi_1}+1/\tau^2}<1 \ \ \ \ \forall \ i=1,...,q.$$
The right-side inequality of \eqref{NC.fishfish} can be written as
\begin{eqnarray}
\label{eq:FI_condition}
\sum_{i=1}^{q}\dfrac{1}{\left(1+(\tau^2 \lambda_i^{\Psi_1})^{-1}\right)^2}  > \dfrac{q}{c}.
\end{eqnarray}
Since ${\displaystyle \sum_{i=1}^{q}\dfrac{1}{\left(1+(\tau^2 \lambda_i^{\Psi_1})^{-1}\right)^2} < q \dfrac{1}{\left( 1+\left( \tau^2  \lambda_{\sup}^{\Psi_1}\right) ^{-1}\right) ^2}}$,
a necessary condition to ensure \eqref{eq:FI_condition} is satisfied is given by
$$ \dfrac{1}{1+\left( \tau^2 \lambda_{\sup}^{\Psi_1}\right) ^{-1} }>\dfrac{1}{\sqrt{c}},$$ noting that $\lambda_{\sup}^{\Psi_1}=\| \Sigma^{-1/2} H  H^T \Sigma^{-1/2}\|_2$. 
\end{proof}

\begin{proof}[Proof of Proposition \ref{prop4bis}] 
Consider that $\tilde{Y}$ follows the general linear form
\begin{eqnarray*}
\tilde{Y} & = & HX + a + U  \ \sim \ \tilde{f}
\end{eqnarray*}
without assumption on $(H,a,S)\in\R^{q,p}\times \R^q \times S^{++}_q(\R)$, apart $S$ being independent on $X$. Then
$\tilde{f}$ is the  ${\cal{N}}(H\mu + a, H\Gamma H^T + S)$ distribution  
and $\arg\min_{f \in {\cal{L}}_{p,q}} \mbox{KL}(f_Y,f)  =  \arg\min_{H,a,S} \mbox{KL}(f_Y, \tilde{f})$.
Since $\tilde{f}$ belongs to the exponential family, it is well known that the KL divergence is convex in the canonic parametrization of $\tilde{f}$, and since it is Gaussian, that this so-called {\it inclusive} KL minimization provides a unique solution for $\tilde{f}$ (see for instance \cite{Minka2001}), by matching  the two first moments of $\tilde{Y}$ and $Y$:
\begin{eqnarray}
a & = & \E[Y]-H\mu \ = \ \E[Y^*] - H\mu \nonumber \\
H\Gamma H^T + S & = & \mbox{Cov}(Y) \ = \ \mbox{Cov}(Y^*) - \Sigma \ \succeq  \ S \succeq 0. \label{constraint.cov.toto}
\end{eqnarray}
Then the result is proved, provided there exists $H\in\R^{q,p}$ such that the covariance constraint (\ref{constraint.cov.toto}) be verified. 
To avoid heavy notations, denote temporarily $Z= \mbox{Cov}(Y)-S$ and assume $Z\succeq 0$. Since $Z$ is symmetric, there exists an orthogonal matrix $O_Z$ and a $q-$diagonal matrix $\Lambda_Z=(\Lambda_{j,j})$ of positive or null terms such that $Z  =  O_Z \Lambda_Z O^T_Z$. Provided $q\leq p$, there always exists $\tilde{\Sigma}_Z\in\R^{q,p}$ such that $\Lambda_Z=\tilde{\Sigma}_Z \tilde{\Sigma}^T_Z$. It is enough to define $\tilde{\Sigma}_Z=(\tilde{\sigma}_{i,j})_{i,j}$, with $1\leq i \leq q$ and $1\leq j \leq p$, such that $\tilde{\sigma}_{j,j}=\sqrt{\Lambda_{j,j}}$, $\tilde{\sigma}_{i,i}=0$ for $i>q$ and $\tilde{\sigma}_{i,j}=0$ for $i\neq j$. Now define in $\R^{q,p}$ 
\begin{eqnarray}
H_{\Gamma,Z} & = & O_Z \tilde{\Sigma}_Z  \Gamma^{-1/2}. \label{sol.H}
\end{eqnarray}
Then, since $\Gamma$ is symmetric, strictly positive definite,  $H_{\Gamma,Z} \Gamma H_{\Gamma,Z}^T  =  O_Z  \tilde{\Sigma}_Z \tilde{\Sigma}^T_Z O^T_Z$ $=
 O_Z \Lambda_Z O^T_Z =  Z$, then  the covariance constraint (\ref{constraint.cov.toto}) is verified. \\
\end{proof} 

\begin{proof}[Proof of Corollary \ref{final.corollary}]
Using the notations of the proof of Proposition \ref{prop4bis}, the matrix $\Psi_{\alpha}$ in previous developments must be replaced by 
 $$
 \Psi_{\Gamma,Z}=(\Sigma+ S) ^{-1/2}H_{\Gamma,Z}H_{\Gamma,Z}^T(\Sigma+S)^{-1/2}.
 $$
 Simplifying the problem to a constraint over a single eigenvalue of $\Gamma$, namely assuming $\Gamma=\tau^2 I_p$, then the adaptation of  (\ref{NC.SC.without.com}), replacing $\Psi_1$ by $\Psi$, is straightforward.
\end{proof}

\section{Computational details}
\label{computational.details}

This section provides additional details on the numerical experiments described in Section \ref{numeric}.

\subsection{Linear surrogate modeling}\label{details}\label{linear.append}

In a real-world setting where the computer code is expensive to evaluate, the inverse problem resolution may become computationally intractable. A worka\-round is to approximate the computer model by a surrogate model much cheaper to run so the inverse problem is solved in a reasonable time. While a considerable literature is already available on Gaussian processes meta-models \cite{Kennedy2001, Fu2015, Fu2016}, the situation considered here is focuses on the simpler  case where the model to invert is approximated by a linear surrogate. As stated in $\S$ \ref{linear.surrg}, a maximin Latin Hypercube Design of $N=2700$ values of $X$ was sampled within the cubic domain ${\cal{A}}$ and the corresponding values $Y^{(k)}=g(X,d^{(k)})$ were computed, given a sampled value of $d^{(k)}$, to fit each linear surrogate.  Such a design maximizes the minimum distance $\delta_D = \min_{i \neq j} \lVert z_i - z_j \rVert$ between data points (see \cite{Koehler1996} and \cite{Damblin2013} for more precisions). 

More precisely, each surrogate was selected from a training generated from $N'=2000$ intervals along each dimension of the domain ${\cal{A}}$, providing $N'$ samples; this is referred to as the \textit{initial design}. At testing, as advocated by \cite{Iooss.Boussouf.ea2010}, $M = 700$ new design points were sequentially generated from the initial design while maintaining the Latin Hypercube initial design property. Note that, unlike \cite{Iooss.Boussouf.ea2010}, the initial design was not filled by choosing the points from a Hammersley sequence that most decreased the design discrepancy. Instead, a new grid was generated, and empty cells corresponding to empty rows and columns were randomly filled until the Latin Hypercube property was respected. Model validation was carried out by computing $RMSE$ and $(CV)RMSE$ values, provided in Table \ref{table:metamodeling}.

\begin{table}[]
    \centering
    \caption{Linear model test scores (averaged on simulations).}
    \label{table:metamodeling}
    \begin{tabular}{ccc}
        \hline
         & $RMSE$ & $(CV)RMSE$ \\
         \hline
         $Y_1$ & 0.15 & 2\% \\
         $Y_2$ & 0.07 & 6\% \\
         \hline
    \end{tabular}
\end{table}

\section{Bayesian stochastic inversion}\label{gelman}

\begin{figure}[hbtp]
    \centering
    \begin{subfigure}{0.45\textwidth}
        \centering
        \includegraphics[width=3.3in]{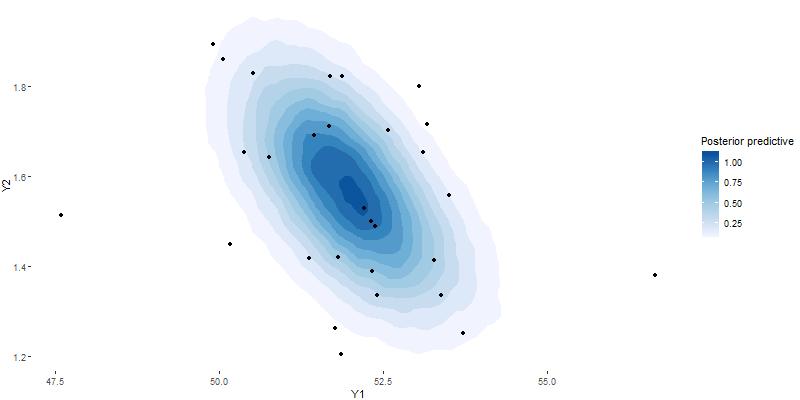}
        \caption{Setting T1.}
        \label{fig:joint-posterior-predictive-Y:a}
    \end{subfigure}
    \begin{subfigure}{0.45\textwidth}
        \centering
        \includegraphics[width=3.3in]{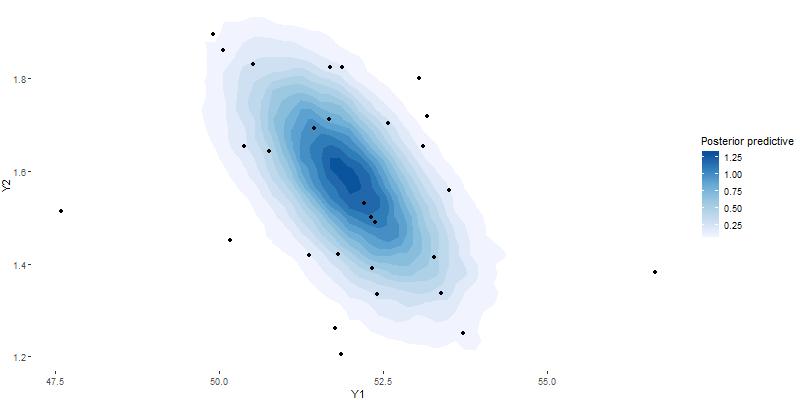}
        \caption{Setting L1.}
        \label{fig:joint-posterior-predictive-Y:b}
    \end{subfigure} \\
    \begin{subfigure}{0.45\textwidth}
        \centering
        \includegraphics[width=3.3in]{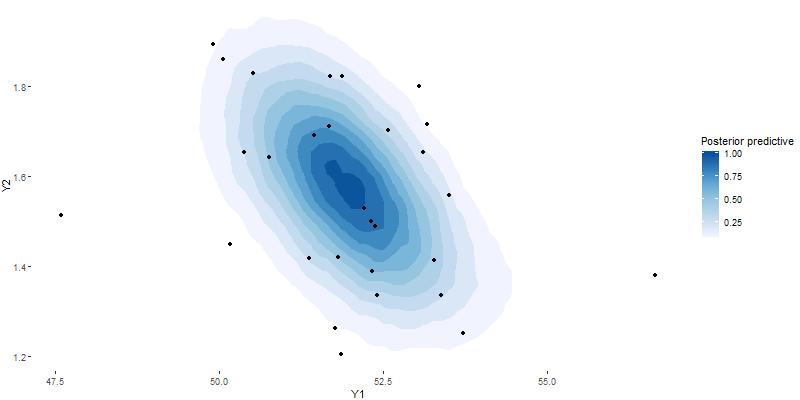}
        \caption{Setting T2.}
        \label{fig:joint-posterior-predictive-Y:c}
    \end{subfigure}
    \begin{subfigure}{0.45\textwidth}
        \centering
        \includegraphics[width=3.3in]{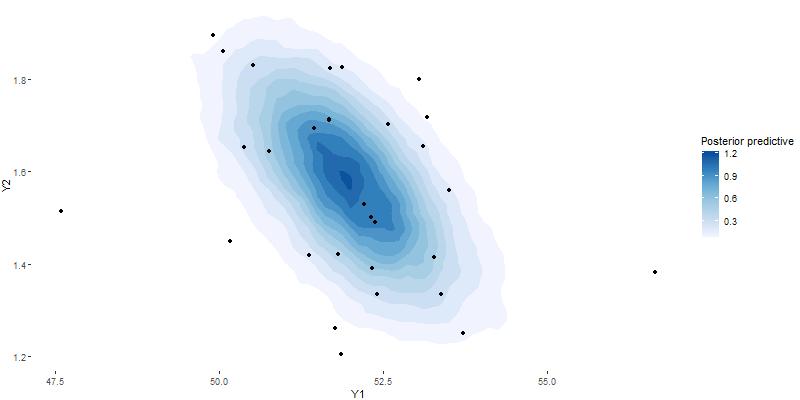}
        \caption{Setting L2.}
        \label{fig:joint-posterior-predictive-Y:d}
    \end{subfigure} \\   
    \hfill\begin{subfigure}{0.5\linewidth}
        \centering
        \includegraphics[width=3.3in]{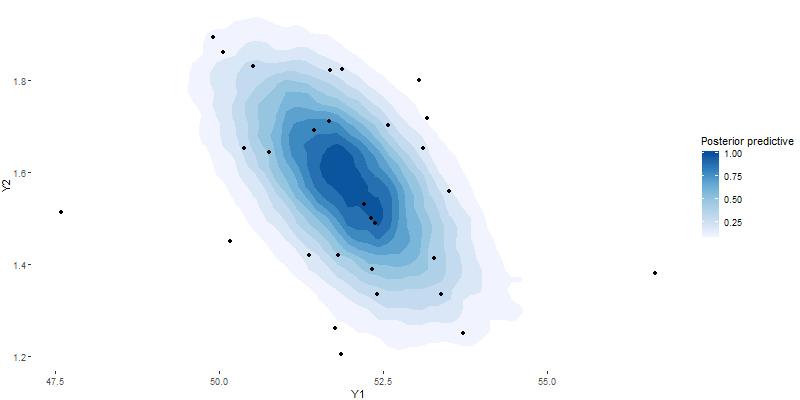}
        \hspace{0.55\textwidth} \caption{Setting L3.}
        \label{fig:joint-posterior-predictive-Y:f}
    \end{subfigure} \\
    \caption{Examples of the joint posterior predictive distribution of $ Y$ obtained under the considered settings. The black dots represents simulated observations $Y$. The left column uses the original model $g$ and the right column uses a linear surrogate.}
    \label{fig:joint-posterior-predictive-Y}
\end{figure}

The Metropolis-Hasting-within-Gibbs algorithm used for the Bayesian computation makes use of three parallel chains, a usual choice that is required to compute convergence diagnostics as Gelman-Rubin $\hat{R}_G$ and Brooks-Gelman statistics  \cite{Gelman1998} . Convergence was firstly monitored by visual inspection in this two-dimensional setting, and the usual rule of thumb $\hat{R}_G < 1.02$ was found to be consistent with the reach of stationarity, after (on averaged) 6000 iterations,  
and by comparing the posterior predictive distributions on $Y$ with simulations (see an illustration on Figures
 \ref{fig:joint-posterior-predictive-Y:a} to \ref{fig:joint-posterior-predictive-Y:f}).

\hspace{-2cm}


\clearpage
\bibliographystyle{plain}
\bibliography{biblio}

\end{document}